\documentclass{amsart}
\usepackage{times,amssymb,amscd,verbatim,graphics,axodraw}
\usepackage[stdtext,vcentermath]{youngtab}
\usepackage[all]{xy}

\numberwithin{equation}{section}

\newtheorem{prop}{Proposition}
\newtheorem{theorem}[prop]{Theorem}
\newtheorem{corollary}[prop]{Corollary}
\newtheorem{lemma}[prop]{Lemma}

\theoremstyle{definition}
\newtheorem{definition}[prop]{Definition}
\newtheorem{example}[prop]{Example}
\newtheorem{remark}[prop]{Remark}

\numberwithin{prop}{section}

\newcommand{\A}{\mathcal{A}}
\newcommand{\cc}{cc}
\newcommand{\Conf}{\mathrm{C}}
\newcommand{\cb}{\overline{c}}
\newcommand{\D}{\mathcal{D}}
\newcommand{\Dp}{\Delta p}
\newcommand{\Dt}{\overleftarrow{D}}
\newcommand{\HH}{\mathcal{H}}
\newcommand{\Hom}{\mathrm{Hom}}
\newcommand{\s}{\overline{s}}
\newcommand{\Jb}{\overline{J}}
\newcommand{\la}{\lambda}
\newcommand{\lab}{\overline{\la}}
\newcommand{\La}{\Lambda}
\newcommand{\lb}{\mathrm{lb}}
\newcommand{\lh}{\mathrm{lh}}
\newcommand{\lm}{\la^-}
\newcommand{\ls}{\mathrm{ls}}
\newcommand{\nub}{\overline{\nu}}
\newcommand{\Path}{\mathcal{P}}
\newcommand{\qbin}[2]{\genfrac{[}{]}{0pt}{}{#1}{#2}}
\newcommand{\RC}{\mathrm{RC}}
\newcommand{\rcls}{\ls_{rc}}
\newcommand{\rclb}{\lb_{rc}}
\newcommand{\rk}{\mathrm{rk}}
\newcommand{\SA}{\mathcal{SA}}
\newcommand{\tb}{\overline{t}}
\newcommand{\word}{\mathrm{word}}
\newcommand{\wt}{\mathrm{wt}}
\newcommand{\Z}{\mathbb{Z}}
\newcommand{\ellb}{\overline{\ell}}
\newcommand{\ft}{\tilde{f}}
\newcommand{\et}{\tilde{e}}
\newcommand{\nut}{\tilde{\nu}}
\newcommand{\Jt}{\tilde{J}}
\newcommand{\ellt}{\tilde{\ell}}

\begin{document}

\title{New fermionic formula  for unrestricted Kostka polynomials}

\author[L.~Deka]{Lipika Deka}
\address{Department of Mathematics, University of California, One Shields
Avenue, Davis, CA 95616-8633, U.S.A.}
\email{deka@math.ucdavis.edu}
\urladdr{http://www.math.ucdavis.edu/\~{}deka}

\author[A.~Schilling]{Anne Schilling}
\email{anne@math.ucdavis.edu}
\urladdr{http://www.math.ucdavis.edu/\~{}anne}

\thanks{\textit{Date:} August 2005}
\thanks{Partially supported by NSF grants DMS-0200774 and DMS-0501101.}

\begin{abstract}
A new fermionic formula for the unrestricted Kostka polynomials of type
$A_{n-1}^{(1)}$  is presented. This formula is different from the one given by 
Hatayama et al. and is valid for all crystal paths based on Kirillov--Reshetihkin
modules, not just for the symmetric and anti-symmetric case. The fermionic formula 
can be interpreted in terms of a new set of unrestricted rigged configurations. For 
the proof a statistics preserving  bijection from this new set of unrestricted rigged
configurations to the set of  unrestricted crystal paths is given which generalizes
a bijection of Kirillov and Reshetikhin.
\end{abstract}

\maketitle

\section{Introduction}
The Kostka numbers $K_{\lambda\mu}$, indexed by the two partitions $\lambda$ and
$\mu$, play an important role in symmetric function theory, representation theory, 
combinatorics, invariant theory and mathematical physics. The Kostka polynomials 
$K_{\lambda\mu}(q)$ are $q$-analogs of the Kostka numbers. There are several
combinatorial definitions of the Kostka polynomials. For example Lascoux and 
Sch{\"u}tzenberger~\cite{LS:1978} proved that the Kostka polynomials are generating
functions of semi-standard tableaux of shape $\lambda$ and content $\mu$
with charge statistic. In~\cite{NY:1997} the Kostka polynomials are expressed
as generating function over highest-weight crystal paths with energy statistics. Crystal
paths are elements in tensor products of finite-dimensional crystals.
Dropping the highest-weight condition yields unrestricted Kostka 
polynomials~\cite{HKOTT:2002,HKKOTY:1999,HKOTY:1999,SW:1999}. In the
$A_1^{(1)}$ setting, unrestricted Kostka polynomials or $q$-supernomial
coefficients were introduced in~\cite{SW:1997} as $q$-analogs of the coefficient of 
$x^a$ in the expansion of $\prod_{j=1}^N (1+x+x^2+\cdots+x^j)^{L_j}$. 
An explicit formula for the $A_{n-1}^{(1)}$ unrestricted Kostka polynomials for 
completely symmetric and completely antisymmetric crystals was proved 
in~\cite{HKKOTY:1999,Kir:2000}. This formula is called fermionic as it is a 
manifestly positive expression.  

In this paper we give a new explicit fermionic formula for the unrestricted Kostka 
polynomials for all Kirillov--Reshetikhin crystals of type $A_{n-1}^{(1)}$.
This fermionic formula can be naturally interpreted in terms of a new set
of unrestricted rigged configurations for type $A_{n-1}^{(1)}$.
Rigged configurations are combinatorial objects originating from the Bethe
Ansatz, that label solutions of the Bethe equations. The simplest version of
rigged configurations appeared in Bethe's original paper~\cite{Bethe:1931}
and was later generalized by Kerov, Kirillov and Reshetikhin~\cite{KKR:1986,KR:1988} 
to models with $\mathrm{GL}(n)$ symmetry.
Since the solutions of the Bethe equations label highest weight vectors, one expects 
a bijection between rigged configurations and semi-standard Young tableaux in 
the $\mathrm{GL}(n)$ case. Such a bijection was given in~\cite{KR:1988,KSS:2002}. 
Here we extend this bijection to a bijection $\Phi$ between the new set of unrestricted
rigged configurations and unrestricted paths. It should be noted that $\Phi$ is 
defined algorithmically. In~\cite{S:2005} the bijection was established in a different 
manner by constructing a crystal structure on the set of rigged configurations.
Here we show that the crystal structures are compatible under the algorithmically
defined $\Phi$ and use this to prove that $\Phi$ preserves the statistics.

Recently, fermionic expressions for generating functions of unrestricted paths
for type $A_1^{(1)}$ have also surfaced in connection with box-ball systems.
Takagi~\cite{T:2004} establishes a bijection between box-ball systems and a new set of 
rigged configurations to prove a fermionic formula for the $q$-binomial coefficient.
His set of rigged configurations coincides with our set in the type $A_1^{(1)}$ case.
There is a generalization of Takagi's bijection to type $A_{n-1}^{(1)}$ 
case~\cite{KOSTY:2005} . Hence this generalization gives a box-ball interpretation
of the unrestricted rigged configurations.

One of the motivations to seek an explicit expression for unrestricted Kostka
polynomials is their appearance in generalizations of the Bailey 
lemma~\cite{Bailey:1949}. Bailey's lemma is a very powerful method to prove 
Rogers--Ramanujan-type identities. In~\cite{SW:1999} a type $A_n$ generalization 
of Bailey's lemma was conjectured which was subsequently
proven in~\cite{W:2002}. A type $A_2$ Bailey chain, which yields an infinite
family of identities, was given in~\cite{ASW:1999}. The new fermionic formulas
of this paper might trigger further progress towards generalizations of the
Bailey lemma.

The bijection $\Phi$ has been implemented as a C++ program~\cite{Deka:2005}
and has been incorporated into the combinatorics package of MuPAD-Combinat by
Francois Descouens~\cite{MuPAD:2005}.

This paper is structured as follows. In Section~\ref{sec:paths} we review crystals
of type $A_{n-1}^{(1)}$, unrestricted paths and the definition of unrestricted Kostka 
polynomials as generating functions of unrestricted paths with energy statistics.
In Section~\ref{sec:RC} we give our new definition of unrestricted rigged configurations
(see Definition~\ref{def:uRC}) and derive from this a fermionic expression for the 
generating function of unrestricted rigged configurations graded by cocharge (see
Section~\ref{sec:fermionic}). The statistic preserving bijection between unrestricted 
paths and unrestricted rigged configurations is established in Section~\ref{sec:main} 
(see Definition~\ref{def:bij} and Theorem~\ref{thm:bij}). As a corolloray this
yields the equality of the unrestricted Kostka polynomials and the fermionic
formula of Section~\ref{sec:RC} (see Corolloray~\ref{cor:X=M}). The result that
the crystal structures on paths and rigged configurations are compatible under
$\Phi$ is stated in Theorem~\ref{thm:commute}. Most of the technical proofs are relegated 
to three appendices. An extended abstract of this paper can be found in~\cite{DS:2005}.

\section{Unrestricted paths and Kostka polynomials} \label{sec:paths}

\subsection{Crystals $B^{r,s}$ of type $A_{n-1}^{(1)}$}
Kashiwara~\cite{Kash:1990}  introduced the notion of crystals and crystal graphs 
as a combinatorial means to study representations of quantum algebras associated 
with any symmetrizable Kac--Moody algebra. In this paper we only consider the 
Kirillov--Reshetikhin crystal $B^{r,s}$ of type $A_{n-1}^{(1)}$ and hence restrict
to this case here.

As a set, the crystal $B^{r,s}$ consists of all column-strict Young tableaux of shape 
$(s^r)$ over the alphabet $\{1,2,\ldots,n\}$. As a crystal associated to the
underlying algebra of finite type $A_{n-1}$, $B^{r,s}$ is isomorphic to the highest
weight crystal with highest weight $(s^r)$. We will define the classical crystal
operators explicitly here. The affine crystal operators $e_0$ and $f_0$ are given
explicitly in~\cite{Sh:2002}. Since we do not use these operators in this
paper we will omit the details.

Let $I=\{1,2,\ldots,n-1\}$ be the index set for the vertices of the Dynkin
diagram of type $A_{n-1}$, $P$ the weight lattice, $\{\La_i\in P \mid i\in I\}$
the fundamental roots, $\{\alpha_i\in P \mid i\in I\}$ the 
simple roots, and $\{h_i\in \Hom_{\Z}(P,\Z) \mid i\in I\}$ the simple coroots.  
As a type $A_{n-1}$ crystal, $B=B^{r,s}$ is equipped with maps 
$e_i,f_i: B\longrightarrow B\cup\{0\}$ and  $\wt: B\longrightarrow P$ 
for all $ i\in I$ satisfying 
\begin{equation*}
\begin{split}
 &f_i(b)=b' \Leftrightarrow e_i(b')=b \text{ if } b,b'\in B\\
 &\wt(f_i(b))=\wt(b)-\alpha_i \text{ if } f_i(b)\in B\\
 &\langle h_i,\wt(b) \rangle=\varphi_i(b)-\varepsilon_i(b),
\end{split}
\end{equation*}
where $\langle \cdot,\cdot \rangle$ is the natural pairing. The maps $f_i,e_i$ are known as
the Kashiwara operators. Here for $b\in B$,
\begin{align*}
 \varepsilon_i(b)&=\max\{k\ge 0 \mid e_i^k(b)\ne 0\}\\
 \varphi_i(b)&=\max\{k\ge 0 \mid f_i^k(b)\ne 0\}.
\end{align*}
Note that for type $A_{n-1}$, $P=\Z^{n}$ and $ \alpha_i=\epsilon_i-\epsilon_{i+1}$
where $\{\epsilon_i \mid i\in I\}$ is the standard basis in $P$. Here $\wt(b)=
(w_1,\ldots,w_n)$ is the weight of $b$ where $w_i$ counts the number of letters 
$i$ in $b$. 

Following~\cite{KN:1994} let us give the action of $e_i$ and $f_i$ 
for $i\in I$.  Let $b \in B^{r,s}$ be a tableau of shape $(s^r)$. The row word of $b$ is 
defined by $\word(b)=w_r\cdots w_2w_1$ where $w_k$ is the word obtained by reading the 
$k$-th row of $b$ from left to right. 
To find $f_i(b)$ and $e_i(b)$ we only consider the subword consisting of the letters $i$ 
and $i+1$ in the word of $b$. First view each $i+1$ in the subword as an opening bracket and
each $i$ as a closing bracket. Then we ignore each adjacent pair of matched 
brackets  successively. At the end of this process we are left with a subword of the
form $i^p(i+1)^q$. If $p>0$ (resp. $q>0$) then $f_i(b)$ (resp. $e_i(b)$) 
is obtained from $b$ by replacing the unmatched subword $i^p(i+1)^q$ 
by $i^{p-1}(i+1)^{q+1}$ (resp. $i^{p+1}(i+1)^{q-1}$). If $p=0$ (resp. $q=0$) then $f_i(b)$ 
(resp. $e_i(b)$) is undefined and we write $f_i(b)=0$ (resp. $e_i(b)=0$). 
 
 A crystal $B$ can be viewed as a directed edge-colored graph whose vertices 
 are the elements of $B$, with a directed edge from $b$
 to $b'$ labeled $i\in I$, if and only if $f_i(b)=b'$. This directed graph is known as
 the crystal graph.   
 \begin{example}
 The crystal graph for $B=B^{1,1}$ is given in Figure~\ref{fig:1}.
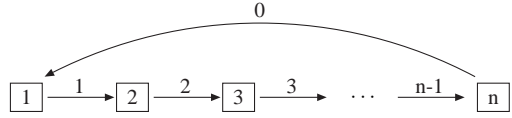
\begin{figure}
\scalebox{.8}{
\begin{picture}(250,62)(-10,-12)
\BText(0,0){1} \LongArrow(10,0)(40,0) \BText(50,0){2}
\LongArrow(60,0)(90,0) \BText(100,0){3} \LongArrow(110,0)(140,0)
\Text(160,0)[]{$\cdots$} \LongArrow(175,0)(205,0)
\BText(220,0){n} \LongArrowArc(110,-181)(216,62,118)
\PText(25,2)(0)[b]{1} \PText(75,2)(0)[b]{2} \PText(125,2)(0)[b]{3}
\PText(190,2)(0)[b]{n-1} \PText(110,38)(0)[b]{0}
\end{picture}}
\caption{Crystal $B^{1,1}$.\label{fig:1}}
\end{figure}
\end{example}
 
Given two crystals $B$ and $B'$, we can also define a new crystal by taking the tensor 
product $B\otimes B'$. As a set $B\otimes B'$ is just the Cartesian product of the
sets $B$ and $B'$. The weight function $\wt$ for $b\otimes b' \in B\otimes B'$
is $\wt(b\otimes b')=\wt(b)+\wt(b')$ and the Kashiwara operators $e_i$, $f_i$ are defined 
as follows
\begin{equation*}
\begin{split}
e_i(b\otimes b')&=
\begin{cases} e_ib\otimes b' & \text{if $\varepsilon_i(b)>\varphi_i(b')$,}\\
              b\otimes e_i b' & \text{otherwise,}\end{cases}\\
f_i(b\otimes b')&=
\begin{cases} f_ib\otimes b' & \text{if $\varepsilon_i(b)\ge \varphi_i(b')$,}\\
              b\otimes f_i b' & \text{otherwise.}\end{cases}
\end{split}
\end{equation*}
This action of $f_i$ and $e_i$ on the tensor product is compatible with the 
previously defined action on $\word(b\otimes b')=\word(b)\word(b')$.

\begin{example}
Let $i=2$ and  
\begin{equation*}
b\;=\; \young(12,23) \otimes \young(23,34,45).
\end{equation*}
Then $\word(b)=2312453423$, the relevant subword is $23-2--3-23$, and the unmatched 
subword is $2--------3$. Hence
\begin{equation*}
f_2(b)\;=\;\young(12,33) \otimes \young(23,34,45)\quad \text{and} \quad
e_2(b)\;=\;\young(12,23) \otimes \young(22,34,45). 
\end{equation*}
\end{example} 

\subsection{Unrestricted paths}
$A_{n-1}^{(1)}$-unrestricted Kostka polynomials or supernomial coefficients were first 
introduced in~\cite{SW:1999} as generating functions of unrestricted paths graded by an 
energy function. An unrestricted path is an element in the tensor product of crystals
$B=B^{r_k,s_k}\otimes B^{r_{k-1},s_{k-1}}\otimes \cdots \otimes B^{r_1,s_1}$.

Let $\la=(\la_1,\la_2,\ldots,\la_n)$ be an $n$-tuple of nonnegative integers. The set 
of \textbf{unrestricted paths} is defined as
\begin{equation*}
\Path(B,\la)=\{b\in B\mid \wt(b)=\la\}.
\end{equation*}

\begin{example}
For $B=B^{1,1}\otimes B^{2,2}\otimes B^{3,1}$ of type $A_3$ and $\la=(2,3,1,2)$ the
path
\begin{equation*}
b\;=\; \young(2) \otimes \young(12,24) \otimes \young(1,3,4)
\end{equation*}
is in $\Path(B,\la)$.
\end{example}  

There exists a crystal isomorphism $R:B^{r,s}\otimes B^{r',s'} \to B^{r',s'} \otimes B^{r,s}$,
called the \textbf{combinatorial $R$-matrix}. Combinatorially it is given as follows. 
Let $b\in B^{r,s}$ and $b'\in B^{r',s'}$.
The product $b\cdot b'$ of two tableaux is defined as the Schensted insertion of $b'$ into
$b$. Then $R(b\otimes b')=\tilde{b}'\otimes \tilde{b}$ is the unique pair of tableaux
such that $b\cdot b'=\tilde{b}'\cdot\tilde{b}$.

The \textbf{local energy function} $H:B^{r,s}\otimes B^{r',s'}\to \Z$ is defined as
follows. For $b\otimes b'\in B^{r,s}\otimes B^{r',s'}$, $H(b\otimes b')$ is the number
of boxes of the shape of $b\cdot b'$ outside the shape obtained by concatenating 
$(s^r)$ and $({s'}^{r'})$.

\begin{example}
For 
\begin{equation*}
b\otimes b'= \young(12,24) \otimes \young(1,3,4)
\end{equation*}
we have 
\begin{equation*}
b\cdot b' = \young(113,224,4) = \young(1,2,4) \cdot \young(13,24) = \tilde{b}'\cdot\tilde{b}.
\end{equation*}
so that
\begin{equation*}
R(b\otimes b')=\tilde{b}'\otimes\tilde{b}=\young(1,2,4) \otimes \young(13,24).
\end{equation*}
Since the concatentation of $\yng(2,2)$ and $\yng(1,1,1)$ is $\yng(3,3,1)$, the local
energy function $H(b\otimes b')=0$.
\end{example}

Now let $B=B^{r_k,s_k}\otimes \cdots\otimes B^{r_1,s_1}$ be a $k$-fold tensor
product of crystals. The \textbf{tail energy function} $\Dt:B\to \Z$ is given by
\begin{equation*}
  \Dt(b) = \sum_{1\le i<j\le k} H_{j-1} R_{j-2} \dotsm R_{i+1} R_i(b),
\end{equation*}
where $H_i$ (resp. $R_i$) is the local energy function (resp. combinatorial $R$-matrix)
acting on the $i$-th and $(i+1)$-th tensor factors of $b\in B$.

\begin{definition}
The $q$-\textbf{supernomial coefficient} or the \textbf{unrestricted Kostka polynomial} is 
the generating function of unrestricted paths graded by the tail energy function
\begin{equation*}
X(B,\la)=\sum_{b\in \Path(B,\la)} q^{\Dt(b)}.
\end{equation*}
\end{definition}

\section{Unrestricted rigged configurations and fermionic formula}
\label{sec:RC}
Rigged configurations are combinatorial objects invented to label
the solutions of the Bethe equations, which give the eigenvalues of
the Hamiltonian of the underlying physical model~\cite{Bethe:1931}.
Motivated by the fact that representation theoretically the eigenvectors
and eigenvalues can also be labelled by Young tableaux, Kirillov
and Reshetikhin~\cite{KR:1988} gave a bijection between tableaux
and rigged configurations. This result and generalizations thereof
were proven in~\cite{KSS:2002}.

In terms of crystal base theory, the bijection is between highest
weight paths and rigged configurations. The new result of this paper
is an extension of this bijection to a bijection between unrestricted paths
and a new set of rigged configurations. The new set of unrestricted rigged
configurations is defined in this section, whereas the bijection is given in 
section~\ref{sec:main}.
In~\cite{S:2005}, a crystal structure on the new set of unrestricted rigged 
configurations is given, which provides a different description of the bijection.

\subsection{Unrestricted rigged configurations}
Let $B=B^{r_k,s_k}\otimes \cdots \otimes B^{r_1,s_1}$ and denote by
$L=(L_i^{(a)}\mid (a,i)\in \HH)$ the multiplicity array of $B$, where
$L_i^{(a)}$ is the multiplicity of $B^{a,i}$ in $B$. Here 
$\HH=I \times \Z_{>0}$ and $I=\{1,2,\ldots,n-1\}$ is the index set of the
Dynkin diagram $A_{n-1}$.
The sequence of partitions $\nu=\{\nu^{(a)}\mid a\in I \}$ is a
\textbf{$(L,\la)$-configuration} if
\begin{equation}\label{eq:size1}
\sum_{(a,i)\in\HH} i m_i^{(a)} \alpha_a = \sum_{(a,i)\in\HH} i
L_i^{(a)} \La_a- \la,
\end{equation}
where $m_i^{(a)}$ is the number of parts of length $i$ in partition
$\nu^{(a)}$. Note that we do not require $\la$ to be a dominant weight here.
The \textbf{(quasi-)vacancy number} of a configuration is defined as
\begin{equation*}
p_i^{(a)}=\sum_{j\ge 1} \min(i,j) L_j^{(a)}
 - \sum_{(b,j)\in \HH} (\alpha_a | \alpha_b) \min(i,j)m_j^{(b)}.
\end{equation*}
Here $(\cdot | \cdot )$ is the normalized invariant form on the weight lattice $P$
such that $(\alpha_i | \alpha_j)$ is the Cartan matrix. Let
$\Conf(L,\la)$ be the set of all $(L,\la)$-configurations.
We call $p_i^{(a)}$ quasi-vacancy number to indicate that they can actually
be negative in our setting. For the rest of the paper we will simply call
them vacancy numbers. 

When the dependence of $m_i^{(a)}$ and $p_i^{(a)}$ on the configuration
$\nu$ is crucial, we also write $m_i^{(a)}(\nu)$ and $p_i^{(a)}(\nu)$,
respectively.

In the usual setting a rigged configuration $(\nu,J)$ consists of a configuration
$\nu\in \Conf(L,\la)$ together with a double sequence of partitions
$J=\{J^{(a,i)}\mid (a,i)\in\HH \}$ such that the partition
$J^{(a,i)}$ is contained in a $m_i^{(a)}\times p_i^{(a)}$ rectangle.
In particular this requires that $p_i^{(a)}\ge 0$. For unrestricted paths
we need a bigger set, where the lower bound on the parts in $J^{(a,i)}$
can be less than zero.

To define the lower bounds we need the following notation. Let 
$\la'=(c_1,c_2,\ldots,c_{n-1})^t$ where $c_k=\la_{k+1}+\la_{k+2}+\cdots+\la_n$.
We also set $c_0=c_1$.
Let $\A(\la')$ be the set of tableaux of shape $\la'$ such that the entries 
in column $k$ are from the set $\{1,2,\ldots,c_{k-1}\}$ and are strictly decreasing along 
each column.

\begin{example} For $n=4$ and $\la=(0,1,1,1)$, the set $\A(\la')$
consists of the following tableaux
\begin{equation*}
\young(332,22,1) \quad \young(332,21,1) \quad \young(322,21,1) \quad 
\young(331,22,1) \quad \young(331,21,1) \quad \young(321,21,1).
\end{equation*}
\end{example}
Note that each $t\in \A(\la')$ is weakly decreasing along each row. This is
due to the fact that $t_{j,k}\ge c_k-j+1$ since column $k$ of height $c_k$ is 
strictly decreasing and $c_k-j+1 \ge t_{j,k+1}$ since the entries in column $k+1$ 
are from the set $\{1,2,\ldots,c_{k}\}$.

Given $t\in\A(\la')$, we define the \textbf{lower bound} as
\begin{equation*}
M_i^{(a)}(t)=-\sum_{j=1}^{c_a} \chi(i\ge t_{j,a})
+\sum_{j=1}^{c_{a+1}} \chi(i\ge t_{j,a+1}),
\end{equation*}
where $t_{j,a}$ denotes the entry in row $j$ and column $a$ of $t$,
and $\chi(S)=1$ if the the statement $S$ is true and $\chi(S)=0$ otherwise.

Let $M,p,m\in \Z$ such that $m\ge 0$.
A $(M,p,m)$-quasipartition $\mu$ is a tuple of integers $\mu=(\mu_1,\mu_2,\ldots,\mu_m)$
such that $M\le \mu_m\le \mu_{m-1}\le \cdots\le \mu_1\le p$. Each $\mu_i$ is called
a part of $\mu$. Note that for $M=0$ this would be a partition with at most $m$ parts each 
not exceeding $p$.

\begin{definition} \label{def:uRC}
An \textbf{unrestricted rigged configuration} $(\nu,J)$ associated to a multiplicity
array $L$ and weight $\la$ is a configuration
$\nu\in\Conf(L,\la)$ together with a sequence $J=\{J^{(a,i)}\mid (a,i)\in\HH\}$
where $J^{(a,i)}$ is a $(M_i^{(a)}(t),p_i^{(a)},m_i^{(a)})$-quasipartition
for some $t\in \A(\la')$. Denote the set of all unrestricted rigged configurations
corresponding to $(L,\la)$ by $\RC(L,\la)$.
\end{definition}

\begin{remark}\mbox{}
\begin{enumerate}
\item
Note that this definition is similar to the definition of level-restricted rigged 
configurations~\cite[Definition 5.5]{SS:2001}. Whereas for level-restricted
rigged configurations the vacancy number had to be modified according to 
tableaux in a certain set, here the lower bounds are modified.
\item 
For type $A_1$ we have $\la=(\la_1,\la_2)$ so that $\A=\{ t\}$
contains just the single tableau
\begin{equation*}
t=\begin{array}{|c|} \hline \la_2\\ \hline \la_2-1\\ \hline \vdots\\ \hline 1\\
\hline \end{array}.
\end{equation*}
In this case $M_i(t)=-\sum_{j=1}^{\la_2} \chi(i\ge t_{j,1})=-i$. This agrees
with the findings of~\cite{T:2004}.
\end{enumerate}
\end{remark}
The quasipartition $J^{(a,i)}$ is called \textbf{singular} if it has a
part of size $p_i^{(a)}$. 
It is often useful to view an (unrestricted) rigged configuration $(\nu,J)$ as a
sequence of partitions $\nu$ where the parts of size $i$ in
$\nu^{(a)}$ are labeled by the parts of $J^{(a,i)}$. The pair
$(i,x)$ where $i$ is a part of $\nu^{(a)}$ and $x$ is a part of
$J^{(a,i)}$ is called a \textbf{string} of the $a$-th rigged
partition $(\nu,J)^{(a)}$. The label $x$ is called a \textbf{rigging}.

\begin{example}
Let $n=4$, $\la=(2,2,1,1)$, $L_1^{(1)}=6$ and all other $L_i^{(a)}=0$. Then
\begin{equation*}
(\nu,J) \;=\; \yngrc(3,-2,1,0) \quad \yngrc(2,0) \quad \yngrc(1,-1)
\end{equation*}
is an unrestricted rigged configuration in $\RC(L,\la)$, where we have written
the parts of $J^{(a,i)}$ next to the parts of length $i$ in partition $\nu^{(a)}$.
To see that the riggings form quasipartitions, let us write
the vacancy numbers $p_i^{(a)}$ next to the parts of length $i$ in partition $\nu^{(a)}$:
\begin{equation*}
\yngrc(3,0,1,3) \quad \yngrc(2,0) \quad \yngrc(1,-1).
\end{equation*}
This shows that the labels are indeed all weakly below the vacancy numbers. For
\begin{equation*}
\young(441,33,2,1) \in \A(\la')
\end{equation*}
we get the lower bounds
\begin{equation*}
\yngrc(3,-2,1,-1) \quad \yngrc(2,0) \quad \yngrc(1,-1),
\end{equation*}
which are less or equal to the riggings in $(\nu,J)$.
\end{example}

Let $B=B^{r_k,s_k}\otimes \cdots \otimes B^{r_1,s_1}$ and $L$
the corresponding multiplicity array. Let $(\nu,J)\in\RC(L,\la)$.
Note that rewritting \eqref{eq:size1} we get
\begin{equation}\label{eq:size}
|\nu^{(a)}| = \sum_{j>a} \la_j -\sum_{j=1}^k s_j \max (r_j-a,0).
\end{equation}
Hence for large $i$, by definition of vacancy numbers we have
\begin{equation}\label{eq:p limit}
\begin{split}
p_i^{(a)}&=|\nu^{(a-1)}|-2|\nu^{(a)}|+|\nu^{(a+1)}|+\sum_j \min(i,j)L_j^{(a)}\\
&=\la_a-\la_{a+1}
\end{split}
\end{equation} 
and
\begin{equation}\label{eq:M limit}
\begin{split}
M_i^{(a)}(t)&=-\sum_{j=1}^{c_a} \chi(i\ge t_{j,a})
+\sum_{j=1}^{c_{a+1}} \chi(i\ge t_{j,a+1})\\
&=-c_a+c_{a+1}=-\la_{a+1}.
\end{split}
\end{equation}
For a given $t\in\A(\la')$ define
\begin{equation*}
\Dp_i^{(a)}(t)=p_i^{(a)}-M_i^{(a)}(t).
\end{equation*}
We write $\Dp_i^{(a)}$ for $\Dp_i^{(a)}(t)$ when there is no cause of confusion.
For large $i$, $\Dp_i^{(a)}(t)=\la_a$.

{}From the definition of $p_i^{(a)}$ one may easily verify that
\begin{equation}\label{eq:p ineq}
-p_{i-1}^{(a)}+2p_i^{(a)}-p_{i+1}^{(a)}\ge m_i^{(a-1)}-2m_i^{(a)}+m_i^{(a+1)}.
\end{equation}
Let $t_{\cdot,a}$ denote the $a$-th column of $t$. Then it follows from the
definition of $M_i^{(a)}(t)$ that
\begin{equation*}
M_i^{(a)}(t)=M_{i-1}^{(a)}(t)-\chi(i\in t_{\cdot,a})+\chi(i\in t_{\cdot,a+1}).
\end{equation*}
Hence \eqref{eq:p ineq} can be rewritten as 
\begin{multline}\label{eq:dp ineq}
-\Dp_{i-1}^{(a)}+2\Dp_i^{(a)}-\Dp_{i+1}^{(a)}
-\chi(i\in t_{\cdot,a})+\chi(i\in t_{\cdot,a+1})\\
+\chi(i+1\in t_{\cdot,a})-\chi(i+1\in t_{\cdot,a+1})
\ge m_i^{(a-1)}-2m_i^{(a)}+m_i^{(a+1)}.
\end{multline}

\begin{lemma}\label{lem:convex}
Suppose that for some $t\in \A(\la')$, $\Dp_i^{(a)}(t)\ge 0$ for all $a\in I$ 
and $i$ such that $m_i^{(a)}>0$. Then there exists a $t'\in\A(\la')$ such that
$\Dp_i^{(a)}(t')\ge 0$ for all $i$ and $a$.
\end{lemma}
\begin{proof}
By definition $\Dp_0^{(a)}(t)=0$ and $\Dp_i^{(a)}(t)=\la_a\ge 0$ for large $i$.
By~\eqref{eq:dp ineq}
\begin{multline}\label{eq:convex}
\Dp_i^{(a)}(t)\ge \frac{1}{2}\bigl\{\Dp_{i-1}^{(a)}(t)+\Dp_{i+1}^{(a)}(t)
+\chi(i\in t_{\cdot,a})-\chi(i\in t_{\cdot,a+1})\\
-\chi(i+1\in t_{\cdot,a})+\chi(i+1\in t_{\cdot,a+1})
+m_i^{(a-1)}+m_i^{(a+1)}\bigr\}
\end{multline}
when $m_i^{(a)}=0$. Hence $\Dp_i^{(a)}(t)<0$ is only possible if
$m_i^{(a-1)}=m_i^{(a+1)}=0$, column $a$ of $t$ contains $i+1$ but no $i$,
and column $a+1$ of $t$ contains $i$ but no $i+1$.
Let $k$ be minimal such that $\Dp_i^{(k)}(t)<0$. Note that $k>1$ since
the first column of $t$ contains all letters $1,2,\ldots, c_1$. 
Let $k'\le k$ be minimal such that $\Dp_i^{(a)}(t)=0$ for all $k'\le a<k$.
Then inductively for $a=k-1,k-2,\ldots,k'$ it follows from~\eqref{eq:convex}
that $m_i^{(a-1)}=0$ and column $a$ of $t$ contains $i+1$ but no $i$. Construct a 
new $t'$ from $t$ by replacing all letters $i+1$ in columns $k',k'+1,\ldots,k$ by $i$.
This accomplishes that $\Dp_j^{(a)}(t')\ge 0$ for all $j$ and $1\le a<k$,
$\Dp_i^{(k)}(t')\ge 0$, and $\Dp_j^{(a)}(t')\ge 0$ for all $a\ge k$ such that
$m_j^{(a)}>0$. Repeating the above construction, if necessary, eventually yields 
a new tableau $t''$ such that finally $\Dp_j^{(a)}(t'')\ge 0$ for all $j$ and $a$.
\end{proof}

\subsection{Fermionic formula} \label{sec:fermionic}
The following statistics can be defined on the set of unrestricted rigged configurations.
For $(\nu,J)\in\RC(L,\la)$ let
\begin{equation*}
\cc(\nu,J)=\cc(\nu)+\sum_{(a,i)\in\HH}|J^{(a,i)}|,
\end{equation*}
where $|J^{(a,i)}|$ is the sum of all parts of the quasipartition $J^{(a,i)}$ and 
\begin{equation*}
\cc(\nu)=\frac{1}{2} \sum_{a,b\in I} \sum_{j,k\ge 1} (\alpha_a | \alpha_b) 
 \min(j,k) m_j^{(a)} m_k^{(b)}.
\end{equation*}

\begin{definition}
The RC polynomial is defined as
\begin{equation*}
M(L,\la)=\sum_{(\nu,J)\in\RC(L,\la)} q^{\cc(\nu,J)}.
\end{equation*}
\end{definition}
The RC polynomial is in fact $S_n$-symmetric in the weight $\la$. This is not obvious
from its definition as both~\eqref{eq:size1} and the lower bounds are not
symmetric with respect to $\la$.

Let $\SA(\la')$ be the set of all nonempty subsets of $\A(\la')$ and set
\begin{equation*}
M_i^{(a)}(S)=\max\{M_i^{(a)}(t) \mid t\in S\} \qquad \text{for $S\in\SA(\la')$.}
\end{equation*}
By inclusion-exclusion the set of all allowed riggings for a given $\nu\in\Conf(L,\la)$ is
\begin{equation*}
\bigcup_{S\in\SA(\la')} (-1)^{|S|+1} \{J\mid \text{$J^{(a,i)}$ is a 
 $(M_i^{(a)}(S),p_i^{(a)},m_i^{(a)})$-quasipartition}\}.
\end{equation*}
The $q$-binomial coefficient $\qbin{m+p}{m}$, defined as
\begin{equation*}
\qbin{m+p}{m}=\frac{(q)_{m+p}}{(q)_m(q)_p}
\end{equation*}
where $(q)_n=(1-q)(1-q^2)\cdots(1-q^n)$, is the generating function of partitions
with at most $m$ parts each not exceeding $p$. Hence the polynomial $M(L,\la)$
may be rewritten as 
\begin{multline*}
M(L,\la)=\sum_{S\in\SA(\la')} (-1)^{|S|+1} \sum_{\nu\in\Conf(L,\la)}
q^{\cc(\nu)+\sum_{(a,i)\in\HH} m_i^{(a)}M_i^{(a)}(S)}\\
\times \prod_{(a,i)\in\HH} \qbin{m_i^{(a)}+p_i^{(a)}-M_i^{(a)}(S)}{m_i^{(a)}}
\end{multline*}
called \textbf{fermionic formula}. This formula is different from the fermionic
formulas of~\cite{HKKOTY:1999,Kir:2000} which exist in the special case when
$L$ is the multiplicity array of $B=B^{1,s_k}\otimes \cdots \otimes B^{1,s_1}$
or $B=B^{r_k,1}\otimes \cdots \otimes B^{r_1,1}$.

\section{Bijection} \label{sec:main}

In this section we define the bijection $\Phi:\Path(B,\la)\to\RC(L,\la)$ from
paths to unrestricted rigged configurations algorithmically. The algorithm generalizes
the bijection of~\cite{KSS:2002} to the unrestricted case. The main result is
summarized in the following theorem.
\begin{theorem}\label{thm:bij} Let $B=B^{r_k,s_k}\otimes \cdots \otimes B^{r_1,s_1}$,
$L$ the corresponding multiplicity array and $\la=(\la_1,\ldots,\la_n)$ a sequence
of nonnegative integers.
There exists a bijection $\Phi:\Path(B,\la)\to\RC(L,\la)$ which preserves the statistics,
that is, $\Dt(b)=\cc(\Phi(b))$ for all $b\in\Path(B,\la)$.
\end{theorem}
A different proof of Theorem~\ref{thm:bij} is given in~\cite{S:2005} by proving 
directly that the crystal structure on rigged configurations and paths coincide.
The results in~\cite{S:2005} hold for all for all simply-laced types, not just type
$A_{n-1}^{(1)}$. Hence Theorem~\ref{thm:bij} holds whenever there is a corresponding 
bijection for the highest weight elements (for example for type $D_n^{(1)}$ for 
symmetric powers~\cite{SS:2005} and antisymmetric powers~\cite{S:2005a}). 
Using virtual crystals and the method of folding Dynkin diagrams, these results can 
be extended to other affine root systems. 
In this paper we use the crystal structure to prove that the statistics is preserved. 
It follows from Theorem~\ref{thm:commute} that the algorithmic definition for $\Phi$
of this paper and the definition of~\cite{S:2005} agree.

An immediate corollary of Theorem~\ref{thm:bij} is the relation between the fermionic 
formula for the RC polynomial of section~\ref{sec:RC} and the unrestricted Kostka 
polynomials of section~\ref{sec:paths}.
\begin{corollary}\label{cor:X=M}
With the same assumptions as in Theorem~\ref{thm:bij}, $X(B,\la)=M(L,\la)$.
\end{corollary}

\subsection{Operations on crystals}
To define $\Phi$ we first need to introduce certain maps on paths and
rigged configurations. These maps correspond to the following operations on crystals:
\begin{enumerate}
\item If $B=B^{1,1}\otimes B'$, let $\lh(B)=B'$. This operation is called \textbf{left-hat}.
\item If $B=B^{r,s}\otimes B'$ with $s\ge 2$, let $\ls(B)=B^{r,1}\otimes B^{r,s-1}\otimes B'$.
This operation is called \textbf{left-split}.
\item If $B=B^{r,1}\otimes B'$ with $r\ge 2$, let $\lb(B)=B^{1,1}\otimes B^{r-1,1}\otimes B'$.
This operation is called \textbf{box-split}.
\end{enumerate}
In analogy we define $\lh(L)$ (resp. $\ls(L)$, $\lb(L)$) to be the multiplicity array of 
$\lh(B)$ (resp. $\ls(B)$, $\lb(B)$), if $L$ is the multiplicity array of $B$.
The corresponding maps on crystal elements are given by:
\begin{enumerate}
\item Let $b=c\otimes b'\in B^{1,1}\otimes B'$. Then $\lh(b)=b'$.
\item Let $b=c\otimes b'\in B^{r,s}\otimes B'$, where $c=c_1c_2\cdots c_s$ and $c_i$
denotes the $i$-th column of $c$. Then $\ls(b)=c_1\otimes c_2\cdots c_s\otimes b'$.
\item Let $b=\begin{array}{|c|} \hline b_1\\ \hline b_2\\ \hline \vdots\\ \hline b_r\\ \hline
\end{array}\otimes b'\in B^{r,1}\otimes B'$, where $b_1<\cdots<b_r$.
Then $\lb(b)=\begin{array}{|c|} \hline b_r\\ \hline \end{array} \otimes 
\begin{array}{|c|} \hline b_1\\ \hline \vdots \\ \hline b_{r-1}\\ \hline \end{array} 
\otimes b'$.
\end{enumerate}

In the next subsection we define the corresponding maps on rigged configurations,
and give the bijection in subsection~\ref{ss:bij}.

\subsection{Operations on rigged configurations}
Suppose $L_1^{(1)}>0$. The main algorithm on rigged configurations as defined
in~\cite{KR:1988,KSS:2002} for admissible rigged configurations can be extended
to our setting. For a tuple of nonnegative integers $\la=(\la_1,\ldots,\la_n)$,
let $\lm$ be the set of all nonnegative tuples $\mu=(\mu_1,\ldots,\mu_n)$ such that
$\la-\mu=\epsilon_r$ for some $1\le r\le n$ where $\epsilon_r$ is the canonical $r$-th unit
vector in $\Z^n$. Define $\delta:\RC(L,\la)\to \bigcup_{\mu\in\lm} \RC(\lh(L),\mu)$
by the following algorithm. Let $(\nu,J)\in\RC(L,\la)$. Set $\ell^{(0)}=1$ and repeat the
following process for $a=1,2,\ldots,n-1$ or until stopped. Find the smallest index 
$i\ge \ell^{(a-1)}$ such that $J^{(a,i)}$ is singular. If no such $i$ exists, set 
$\rk(\nu,J)=a$ and stop. Otherwise set $\ell^{(a)}=i$ and continue with $a+1$.
Set all undefined $\ell^{(a)}$ to $\infty$.

The new rigged configuration $(\tilde{\nu},\tilde{J})=\delta(\nu,J)$ is obtained by
removing a box from the selected strings and making the new strings singular
again. Explicitly
\begin{equation*}
 m_i^{(a)}(\tilde{\nu})=m_i^{(a)}(\nu)+\begin{cases}
 1 & \text{if $i=\ell^{(a)}-1$}\\
 -1 & \text{if $i=\ell^{(a)}$}\\
 0 & \text{otherwise.} \end{cases}
\end{equation*}
The partition $\tilde{J}^{(a,i)}$ is obtained from $J^{(a,i)}$ by removing
a part of size $p_i^{(a)}(\nu)$ for $i=\ell^{(a)}$,
adding a part of size $p_i^{(a)}(\tilde{\nu})$ for $i=\ell^{(a)}-1$, 
and leaving it unchanged otherwise. Then $\delta(\nu,J)\in \RC(\lh(L),\mu)$
where $\mu=\la-\epsilon_{\rk(\nu,J)}$.

\begin{prop} \label{prop:delta}
$\delta$ is well-defined.
\end{prop}

The proof is given in Appendix~\ref{appn:delta}.

\begin{example}\label{ex:delta}
Let $L$ be the multiplicity array of $B=B^{1,1}\otimes B^{2,1}\otimes B^{2,3}$
and $\la=(2,2,2,1,1,1)$. Then
\begin{equation*}
(\nu,J)= \yngrc(2,-1,1,0) \quad \yngrc(3,0,1,-1,1,-1) \quad \yngrc(3,0)
\quad \yngrc(2,-1) \quad \yngrc(1,-1) \in \RC(L,\la).
\end{equation*}
Writing the vacancy numbers next to each part instead of the riggings we get
\begin{equation*}
\yngrc(2,-1,1,0) \quad \yngrc(3,0,1,-1,1,-1) \quad \yngrc(3,1)
\quad \yngrc(2,-1) \quad \yngrc(1,-1).
\end{equation*}
Hence $\ell^{(1)}=\ell^{(2)}=1$ and all other $\ell^{(a)}=\infty$, so that
\begin{equation*}
\delta(\nu,J)= \yngrc(2,-1) \quad \yngrc(3,0,1,-1) \quad \yngrc(3,0)
\quad \yngrc(2,-1) \quad \yngrc(1,-1).
\end{equation*}
Also $\cc(\nu,J)=2$.
\end{example}

The inverse algorithm of $\delta$ denoted by $\delta^{-1}$ is defined as follows.
Let $L_1^{(1)}=\overline{L}_1^{(1)}+1, L_i^{(k)}=\overline{L}_i^{(k)}$ for all $i,k\ne1$. 
Let $\lab$ be a weight and $\la=\lab+\epsilon_r$   for some $1\le r\le n$. 
Define $\delta^{-1}: \RC(\overline{L},\lab)\to\RC(L,\la) $ 
by the following algorithm. Let $(\nub,\Jb) \in\RC(\overline{L},\lab)$. 
Let $s^{(r)}=\infty$. For $k=r-1$ down to $1$, select the longest singular string in 
$(\nub,\Jb)^{(k)}$ of length $s^{(k)}$ (possibly of zero length) such that $s^{(k)} \le 
s^{(k+1)}$. With the convention $s^{(0)}=0$ we have $s^{(0)}\le s^{(1)}$ as well. 
$\delta^{-1}(\nub,\Jb)=(\nu,J)$ is obtained from $(\nub,\Jb)$ by adding a box 
to each of the selected strings, and resetting their labels to make them 
singular with respect to the new vacancy number for $\RC(L,\la)$, and 
leaving all other strings unchanged.

\begin{prop}\label{prop:inv delta}
$\delta^{-1}$ is well defined.
\end{prop}
This proposition will also be proved in Appendix~\ref{appn:delta}.

Let $s\ge2$. Suppose $B=B^{r,s}\otimes B'$ and $L$ the corresponding 
multiplicity array. Note that $\Conf(L,\la)\subset \Conf(\ls(L),\la)$. Under this
inclusion map, the vacancy number $p_i^{(a)}$ for $\nu$ increases by
$\delta_{a,r} \chi(i<s)$. Hence there is a well-defined injective map
$\rcls:\RC(L,\la)\rightarrow \RC(\ls(L),\la)$ given by the identity map 
$\rcls(\nu,J)=(\nu,J)$.

Suppose $r\ge2$ and $B=B^{r,1}\otimes B'$ with multiplicity array $L$.
Then there is an injection $\rclb:\RC(L,\la)\to \RC(\lb(L),\la)$ defined by adding 
singular strings of length $1$ to $(\nu,J)^{(a)}$ for $1\le a < r$. Note that the
vacancy numbers remain unchanged under $\rclb$.

\subsection{Bijection}\label{ss:bij}
The map $\Phi:\Path(B,\la)\to\RC(L,\la)$ is defined recursively by various commutative 
diagrams. Note that it is possible to go from $B=B^{r_k,s_k}\otimes 
B^{r_{k-1},s_{k-1}}\otimes \cdots \otimes B^{r_1,s_1}$ to the empty crystal 
via successive application of $\lh$, $\ls$ and
$\lb$.

\begin{definition} \label{def:bij}
Define that map $\Phi:\Path(B,\la)\rightarrow \RC(L,\la)$ such that 
the empty path maps to the empty rigged configuration and such that the following
conditions hold:
\begin{enumerate}
\item \label{bij:1} Suppose $B=B^{1,1} \otimes B'$. Then the following diagram commutes:
\begin{equation*}
\begin{CD}
\Path(B,\la) @>{\Phi}>> \RC(L,\la) \\
@V{\lh}VV @VV{\delta}V \\
\displaystyle{\bigcup_{\mu\in\lm} \Path(\lh(B),\mu)} @>>{\Phi}> 
\displaystyle{\bigcup_{\mu\in\lm}
\RC(\lh(L),\mu)}
\end{CD}
\end{equation*}
\item \label{bij:2} Suppose $B=B^{r,s} \otimes B'$ with $s\ge 2$. Then the 
following diagram commutes:
\begin{equation*}
\begin{CD}
\Path(B,\la) @>{\Phi}>> \RC(L,\la) \\
@V{\ls}VV @VV{\rcls}V \\
\Path(\ls(B),\la) @>>{\Phi}> \RC(\ls(L),\la)
\end{CD}
\end{equation*}
\item \label{bij:3} Suppose $B=B^{r,1} \otimes B'$ with $r\ge2$. Then the 
following diagram commutes:
\begin{equation*}
\begin{CD}
\Path(B,\la) @>{\Phi}>> \RC(L,\la) \\
@V{\lb}VV @VV{\rclb}V \\
\Path(\lb(B),\la) @>>{\Phi}> \RC(\lb(L),\la)
\end{CD}
\end{equation*}
\end{enumerate}
\end{definition}

\begin{prop}\label{prop:bij} The map $\Phi$ of Definition~\ref{def:bij}
is a well-defined bijection.
\end{prop}
The proof is given in Appendix~\ref{appn:phi}.

\begin{example} \label{ex:b}
Let $B=B^{1,1}\otimes B^{2,1}\otimes B^{2,3}$ and $\la=(2,2,2,1,1,1)$. Then
\begin{equation*}
b=\young(3) \otimes \young(1,2) \otimes \young(123,456)\in\Path(B,\la)
\end{equation*}
and $\Phi(b)$ is the rigged configuration $(\nu,J)$ of Example~\ref{ex:delta}.
We have $\Dt(b)=\cc(\nu,J)=2$.
\end{example}

\begin{example}
Let $n=4$, $B=B^{2,2}\otimes B^{2,1}$ and $\la=(2,2,1,1)$. Then the multiplicity
array is $L_1^{(2)}=1,L_2^{(2)}=1$ and $L_i^{(a)}=0$ for all other $(a,i)$. There are 
7 possible unrestricted paths in $\Path(B,\la)$. For each path $b\in\Path(B,\la)$ the 
corresponding rigged configuration $(\nu,J)=\Phi(b)$ together with the tail energy 
and cocharge is summarized below.
\begin{equation*}
\begin{array}{llllll}
b\;=\; \young(11,22) \otimes \young(3,4)
&
\quad (\nu,J) \;=
& \yngrc(1,0) & \yngrc(1,-1,1,-1) & \yngrc(1,0)
&
\Dt(b) =0= \cc(\nu,J)\\[5mm]
b\;=\; \young(11,24) \otimes \young(2,3)
&
\quad (\nu,J) \;=
& \yngrc(1,-1) & \yngrc(1,0,1,0) & \yngrc(1,0)
&
\Dt(b) =1= \cc(\nu,J)\\[5mm]
b\;=\; \young(12,23) \otimes \young(1,4)
&
\quad (\nu,J) \;=
& \yngrc(1,0) & \yngrc(1,0,1,0) & \yngrc(1,-1)
&
\Dt(b) =1= \cc(\nu,J)\\[5mm]
b\;=\; \young(12,24) \otimes \young(1,3)
&
\quad (\nu,J) \;=
& \yngrc(1,0) & \yngrc(1,0,1,-1) & \yngrc(1,0)
&
\Dt(b)=1= \cc(\nu,J)\\[5mm]
b\;=\; \young(13,24) \otimes \young(1,2)
&
\quad (\nu,J) \;=
& \yngrc(1,0) & \yngrc(1,0,1,0) & \yngrc(1,0)
&
\Dt(b)=2= \cc(\nu,J)\\[5mm]
b\;=\; \young(11,23) \otimes \young(2,4)
&
\quad (\nu,J) \;=
& \yngrc(1,-1) & \yngrc(2,0) & \yngrc(1,-1)
&
\Dt(b)=0= \cc(\nu,J)\\[5mm]
b\;=\; \young(12,34) \otimes \young(1,2)
&
\quad (\nu,J) \;=
& \yngrc(1,-1) & \yngrc(2,1) & \yngrc(1,-1)
&
\Dt(b)=1= \cc(\nu,J)
\end{array}
\end{equation*}
The unrestricted Kostka polynomial in this case is $M(L,\la)=2+4q+q^2=X(B,\la)$.
\end{example}

\subsection{Crystal operators on unrestricted rigged configurations}
Let $B=B^{r_k,s_s}\otimes \cdots \otimes B^{r_1,s_1}$ and $L$ be the multiplicity 
array of $B$. Let $\Path(B)=\bigcup_{\la}\Path(B,\la)$ and $\RC(L)=\bigcup_{\la}\RC(L,\la)$. 
Note that the bijection $\Phi$ of Definition~\ref{def:bij} extends to 
a bijection from $\Path(B)$ to $\RC(L)$. Let $f_a$ and $e_a$ for 
$1\le a<n$ be the crystal operators acting on the paths in $\Path(B)$. 
In~\cite{S:2005} analogous operators $\ft_a$ and $\et_a$ for $1\le a<n$ acting on rigged 
configurations in $\RC(L)$ were defined.   

\begin{definition}\cite[Definition 3.3]{S:2005}
\begin{enumerate}
\item
Define $\et_a(\nu,J)$ by removing a box from a string of length $k$ in
$(\nu,J)^{(a)}$ leaving all colabels fixed and increasing the new
label by one. Here $k$ is the length of the string with the smallest
negative rigging of smallest length. If no such string exists,
$\et_a(\nu,J)$ is undefined.
\item
Define $\ft_a(\nu,J)$ by adding a box to a string of length $k$ in
$(\nu,J)^{(a)}$ leaving all colabels fixed and decreasing the new
label by one. Here $k$ is the length of the string with the smallest
nonpositive rigging of largest length. If no such string exists,
add a new string of length one and label -1.
If the result is not a valid unrestricted rigged configuration
$\ft_a(\nu,J)$ is undefined.
\end{enumerate}
\end{definition}

\begin{example}
Let $L$ be the multiplicity array of $B=B^{1,3}\otimes B^{3,2}\otimes B^{2,1}$ 
and let
\begin{equation*}
(\nu,J)= \yngrc(4,-3,1,-1) \quad \yngrc(3,0,1,1) \quad \yngrc(2,-1,1,-1) \in \RC(L).
\end{equation*} 
Then 
\begin{equation*}
\begin{split}
\ft_3(\nu,J)&= \yngrc(4,-3,1,-1) \quad \yngrc(3,1,1,1) \quad \yngrc(3,-2,1,-1)\\
\text{and} \qquad
\et_3(\nu,J)&= \yngrc(4,-3,1,-1) \quad \yngrc(3,-1,1,0) \quad \yngrc(2,1).
\end{split}
\end{equation*} 
\end{example}

Define $\widetilde{\varphi}_a(\nu,J)=\max\{k\ge 0 \mid \ft_a(\nu,J)\neq 0\}$ and
$\widetilde{\varepsilon}_a(\nu,J)=\max\{k\ge 0 \mid \et_a(\nu,J)\neq 0\}$.
The following Lemma is proven in~\cite{S:2005}.
\begin{lemma} \cite[Lemma 3.6]{S:2005}\label{lem:varphi}
Let $(\nu,J)\in \RC(L)$. For fixed $a\in \{1,2,\ldots,n-1\}$, let $p=p_i^{(a)}$ be the 
vacancy number for large $i$ and let $s\le 0$ be the smallest nonpositive label
in $(\nu,J)^{(a)}$; if no such label exists set $s=0$. 
Then $\widetilde{\varphi}_a(\nu,J)=p-s$.
\end{lemma}

\begin{theorem} \label{thm:commute}
Let $B=B^{r_k,s_k}\otimes \cdots \otimes B^{r_1,s_1}$ and $L$ the multiplicity array 
of $B$. Then  the following diagrams commute:
\begin{equation}\label{eq:commute}
\begin{CD}
\Path(B) @>{\Phi}>> \RC(L) \\
@V{f_a}VV @VV{\ft_a}V \\
\Path(B) @>>{\Phi}> \RC(L)
\end{CD}
\qquad \qquad
\begin{CD}
\Path(B) @>{\Phi}>> \RC(L) \\
@V{e_a}VV @VV{\et_a}V \\
\Path(B) @>>{\Phi}> \RC(L).
\end{CD}
\end{equation} 
\end{theorem}
The proof of  Theorem~\ref{thm:commute} is given in Appendix~\ref{appn:crystal}.  
Note that Proposition~\ref{prop:bij} and Theorem~\ref{thm:commute} imply
that the operators $\ft_a,\et_a$ give a crystal structure on $\RC(L)$. 
In~\cite{S:2005} it is shown directly that $\ft_a$ and $\et_a$ define a crystal 
structure on $\RC(L)$.

\subsection{Proof of Theorem~\ref{thm:bij}}
By Proposition~\ref{prop:bij} $\Phi$ is a bijection which proves the first part of
Theorem~\ref{thm:bij}. By Theorem~\ref{thm:commute} the operators $\ft_a$ and $\et_a$ give 
a crystal structure on $\RC(L)$ induced by the crystal structure on $\Path(B)$
under $\Phi$. The highest weight elements are given by the usual rigged configurations 
and highest weight paths, respectively, for which Theorem~\ref{thm:bij} is known to 
hold by~\cite{KSS:2002}. The energy function $\Dt$ is constant on classical
components. By~\cite[Theorem 3.9]{S:2005} the statistics $\cc$ on rigged configurations
is also constant on classical components. Hence $\Phi$ preserves the statistic.

\subsection{Implementation} The bijection $\Phi$ and its inverse have been implemented
as a C++ program. The code is available in~\cite{Deka:2005}. In early stages
of this project these programs have been invaluable to produce data and 
check conjectures regarding the unrestricted rigged configurations.
The progams have also been incorporated into MuPAD-Combinat as a dynamic
module by Francois Descouens~\cite{MuPAD:2005}.
For example, the command
\begin{multline*}
\text{{\tt riggedConfigurations::RcPathsEnergy::}}\\
\text{{\tt fromOnePath([[[3]],[[2],[1]],[[4,5,6],[1,2,3]]])}}
\end{multline*}
calculates $\Phi(b)$ with $b$ as in Example~\ref{ex:b}.

\appendix
\section{Proof of Propositions~\ref{prop:delta} and~\ref{prop:inv delta}}\label{appn:delta}

In this section we prove Propositions~\ref{prop:delta} and~\ref{prop:inv delta},
namely that $\delta$ is a well-defined bijection.
The following remark will be useful.

\begin{remark}\label{remark:new tab}
Let $(\nu,J)$ be admissible with respect to $t\in\A(\la')$.
Suppose that $\Dp_{i-1}^{(k)}(t)+\Dp_{i+1}^{(k)}(t)\ge 1$ and $\Dp_i^{(k)}(t)=m_i^{(k)}(\nu)=0$. 
Then by \eqref{eq:dp ineq} there are five choices for the letters $i$ and $i+1$
in columns $k$ and $k+1$ of $t$:
\begin{enumerate}
\item \label{l:1} $i+1$ in column $k$;
\item \label{l:2} $i+1$ in column $k$ and $k+1$, $i$ in column $k+1$;
\item \label{l:3} $i$ in column $k+1$;
\item \label{l:4} $i$ in column $k$ and $k+1$, $i+1$ in column $k$;
\item \label{l:5} $i+1$ in column $k$, $i$ in column $k+1$.
\end{enumerate}
In cases \ref{l:1} and \ref{l:2} we have $m_i^{(k-1)}(\nu)=0$. Changing letter $i+1$
to $i$ in column $k$ to form a new tableau $t'$ has the effect
$M_i^{(k)}(t')=M_i^{(k)}(t)-1$, $M_i^{(k-1)}(t')=M_i^{(k-1)}(t)+1$ and all
other lower bounds remain unchanged.
In cases \ref{l:3} and \ref{l:4} we have $m_i^{(k+1)}(\nu)=0$. Changing letter $i$
to $i+1$ in column $k+1$ to form a new tableau $t'$ has the effect
$M_i^{(k)}(t')=M_i^{(k)}(t)-1$, $M_i^{(k+1)}(t')=M_i^{(k+1)}(t)+1$ and all
other lower bounds remain unchanged.
Finally in case \ref{l:5} either $m_i^{(k-1)}(\nu)=0$ or $m_i^{(k+1)}(\nu)=0$.
Changing $i+1$ to $i$ in column $k$ (resp. $i$ to $i+1$ in column $k+1$)
has the same effect as in case \ref{l:1} (resp. case \ref{l:3}).

This shows that under the replacement $t\mapsto t'$ we have $\Dp_i^{(k)}(t')>0$
and by Lemma~\ref{lem:convex} $(\nu,J)$ is admissible with respect to some
tableau $t''$.
\end{remark}

Let $\la $ be a weight such that  $\la_r>0$ for a given $1\le r\le n$.
Set $\lab=\la-\epsilon_r$. Recall that $c_k=\la_{k+1}+\la_{k+2}+\cdots+\la_n$ is the 
height of the $k$-th column of $t \in \A(\la')$. Let us define the map 
$\D_r:\A(\la')\to\A(\lab')$ with $\tb=\D_r(t)$ as follows. If $t_{1,r}<c_{r-1}$ then 
\begin{equation}\label{eq:t bar1}
\tb_{i,k} = \begin{cases}
 t_{i+1,k} & \text{for $1\le k\le r-1$ and $1\le i<c_k$,}\\
 t_{i,k} & \text{for $r\le k\le n$ and $1\le i \le c_k$.}
\end{cases}
\end{equation} 
If $t_{1,r}=c_{r-1}$ then there exists $1\le j\le c_r$ such that $t_{i,r} =t_{i-1,r}-1$ for 
$2\le i\le j$ and $t_{j+1,r} < t_{j,r}-1$ if $j<c_r$. In this case
\begin{equation}\label{eq:t bar2}
\tb_{i,k} = \begin{cases}
 t_{i+1,k} & \text{for $1\le k\le r-1$ and $1\le i<c_k$,}\\
 t_{i,r}-1 & \text{for $k=r$ and $1\le i\le j$,}\\
 t_{i,r} &\text{for $k=r$ and $j<i\le c_r$,}\\
 t_{i,k} & \text{for $r<k\le n$ and $1\le i \le c_k$.}
\end{cases}
\end{equation}

Note that by definition the entries of $\D_r(t)$ are strictly decreasing along columns.
Let $\cb_k=\lab_{k+1}+\cdots +\lab_n$. Then we have $\cb_k=c_k-1$ for $1\le k\le r-1$ 
and $\cb_k=c_k$ for $r\le k\le n$. Again by definition
$\tb_{j,1} \in \{1,2,\cdots , \cb_{1}\}$ for all $1\le j\le \cb_1$ and  
$\tb_{j,k} \in \{1,2,\cdots , \cb_{k-1}\}$ for all $2\le j\le \cb_k$ and $1\le k \le n$. 
Therefore,  $\D_r(t)\in\A(\lab')$.

\begin{example}
Let $t=\young(332,21,1)$ and $r=3$. Then $\D_r(t)=\young(211,1)$.
\end{example}

We will use the following lemma and remark in the proofs.
\begin{lemma} \label{lem:biggest part}
Let $B=B^{r_l,s_l}\otimes \cdots \otimes B^{r_1,s_1}$ with 
 $r_l=1=s_l$. Let  $(\overline{\nu},\overline{J})=\delta(\nu,J)$ and let
$\rk(\nu,J)=r$. For $1<k<r$ let $i=t_{1,k}$. Then one of the following 
conditions hold:
\begin{enumerate}
\item $ m_i^{(k)}(\nu)=0$ or
\item $m_i^{(k)}(\nu)=1$, in which case $\delta$ selects the part of length $i$ 
in $\nu^{(k)}$.
\end{enumerate}
\end{lemma} 
\begin{proof} 
Note that $i=t_{1,k} \ge c_k$. 
By \eqref{eq:size} we have $|\nu^{(k)}|\le c_{k}$, so
that either $m_i^{(k)}(\nu)=0$ or $i=c_k$ and
$\nu^{(k)}$ consists of just one part of size $i$. In this case
$m_i^{(k)}(\nu)=1$ and $\delta$ has to select this single part.
\end{proof} 

\begin{remark}\label{remark:rigsize} By \eqref{eq:size} we have
\begin{equation*}
\begin{split}
|\nu^{(r)}|&=|\nu^{(r-1)}|-\la_r+ \sum_{i\ge 1}s_i\chi(r_i\ge r)\\ 
|\nu^{(r+1)}| &=|\nu^{(r-1)}|-\la_r-\la_{r+1}+2\sum_{i\ge 1}s_i\chi(r_i\ge r)-
\sum_{i\ge 1}s_i\delta_{r_i,r}.
\end{split}
\end{equation*}
Note that for $a>0$ 
$$\sum_{i\ge 1}\min(a,i)L_i^{(r)}=\sum_{i\ge 1}s_i\chi(s_i \le a)\delta_{r_i,r}+\sum_
{i\ge 1}a\chi(s_i>a)\delta_{r_i,r}.$$
Then if $|\nu^{(r-1)}|=c_{r-1}-k$ for some $k\ge 0$ it follows that
\begin{equation*}
 -2|\nu^{(r)}|+|\nu^{(r+1)}|+\sum_{i\ge 1}\min(a,i)L_i^{(r)}=-2\la_{r+1}-c_{r+1}+
k-\sum_{i\ge 1}\max(s_i-a,0)\delta_{r_i,r}.
\end{equation*}
\end{remark}

\begin{proof}[Proof of Proposition~\ref{prop:delta}]
To prove that $\delta$ is well-defined it needs to be shown that
$(\nub,\Jb)=\delta(\nu,J)\in\RC(\overline{L},\lab)$. Here 
$\overline{L}$ is given by $\overline{L}_1^{(1)}=L_1^{(1)}-1$,
$\overline{L}_i^{(a)}=L_i^{(a)}$ for all other $i,a$, and 
$\lab=\la-\epsilon_r$ where $r=\rk(\nu,J)$. 

Let us first show that $\lab$ indeed has nonnegative entries. Assume the 
contrary that  $\lab_r<0$. This can happen only if $\la_r=0$ . Suppose 
$t \in \A(\la')$ is such that $M_j^{(k)}(t)\le p_j^{(k)}(\nu)$ for all $j,k$. 
By~\eqref{eq:p limit}, 
$p_i^{(r)}(\nu)=-\la_{r+1}$ for large $i$. Let $\ell$ be the size of the largest 
part in $\nu^{(r)}$, so that $m_j^{(r)}(\nu)=0$ for $j> \ell$. 
By definition of vacancy numbers, $p_i^{(r)}(\nu)\ge p_j^{(r)}(\nu)$ for $i\ge j \ge \ell$. 
Also we have $M_j^{(r)}(t)\ge -\la_{r+1}$ for all $j$. Hence, 
$-\la_{r+1}\le M_{j}^{(r)}(t)\le p_{j}^{(r)}(\nu)\le p_i^{(r)}(\nu)=-\la_{r+1}$ implies 
\begin{equation} \label{eq:peq}
M_i^{(r)}(t)=M_j^{(r)}(t)= p_{j}^{(r)}(\nu)= p_i^{(r)}(\nu) \quad \text{for all $\ell \le j \le i$.}
\end{equation}
This means that the string of length $\ell$ in $(\nu,J)^{(r)}$ is singular and 
$\Dp_j^{(r)}(t)=0$ for all $j\ge \ell$. We claim that $m_j^{(r-1)}(\nu)=0$ for $ j >\ell$. 
By \eqref{eq:dp ineq} we get 
\begin{equation*}
\begin{split}
S:=&-\chi(j\in t_{\cdot,r})+\chi(j\in t_{\cdot,r+1})
+\chi(j+1\in t_{\cdot,r})-\chi(j+1\in t_{\cdot,r+1})\\
\ge & m_j^{(r-1)}(\nu)+m_j^{(r+1)}(\nu)
\end{split}
\end{equation*}
for $j>\ell$. Clearly, $m_j^{(r-1)}(\nu)=0$ unless $1\le S \le 2$. If $S=2$ we have 
$j+1\in t_{\cdot,r}$ and $j\in t_{\cdot,r+1}$ which implies 
$M_j^{(r)}(t)=M_{j+1}^{(r)}(t)+1$, a contradiction to~\eqref{eq:peq}. 
Hence $S=2$ is not possible. Similarly, we can show that $S=1$ is not possible. 
This proves that $m_j^{(r-1)}(\nu)=0$ for $j>\ell$.
Hence $\ell^{(r-1)}\le \ell$ which contradicts the assumption that $r=\rk(\nu,J)$ 
since $(\nu,J)^{(r)}$ has a singular string of length $\ell$. Therefore $\la_r>0$.

Next we need to show that $(\nub,\Jb)$ is admissible, which means that the
parts of $\Jb$ lie between the corresponding lower bound for some $\tb\in \A(\lab')$
and the vacancy number. Let $t\in\A(\la')$ be such that $(\nu,J)$ is admissible
with respect to $t$.
By the same arguments as in the proof of Proposition 3.12 of~\cite{KSS:2002}
the only problematic case is when
\begin{equation}\label{eq:problem}
m_{\ell-1}^{(k)}(\nu)=0, \quad \Dp_{\ell-1}^{(k)}(t)=0, \quad \ell^{(k-1)}<\ell
\quad \text{and $\ell$ finite}
\end{equation}
where $\ell=\ell^{(k)}$. 

Assume that $\Dp_{\ell-2}^{(k)}(t)+\Dp_{\ell}^{(k)}(t)\ge 1$ and \eqref{eq:problem}
holds. By Remark~\ref{remark:new tab} with $i=\ell-1$, there exists a new tableau
$t'$ such that $\Dp_{\ell-1}^{(k)}(t')>0$ so that the problematic case is avoided.

Hence assume that $\Dp_{\ell-2}^{(k)}(t)+\Dp_{\ell}^{(k)}(t)=0$ and \eqref{eq:problem}
holds. Let $\ell'<\ell$ be maximal such that $m_{\ell'}^{(k)}(\nu)>0$. If no such
$\ell'$ exists, set $\ell'=0$. 

Suppose that there exists $\ell'<j<\ell$ such that $\Dp_{j-1}^{(k)}(t)>0$.
Let $i$ be the maximal such $j$. Then by Remark~\ref{remark:new tab}
we can find a new tableau $t'$ such that $\Dp_i^{(k)}(t')>0$ and 
$(\nu,J)$ is admissible with respect to $t'$.
Repeating the argument we can achieve $\Dp_{\ell-1}^{(k)}(t'')>0$ for
some new tableau $t''$, so that the problematic case does not occur.

Hence we are left to consider the case $\Dp_i^{(k)}(t)=0$ for all
$\ell'\le i\le \ell$. If $m_i^{(k-1)}(\nu)=0$ for all $\ell'<i<\ell$, then by the 
same arguments as in the proof of Proposition 3.12 of~\cite{KSS:2002} we arrive 
at a contradition since $\ell^{(k-1)}\le \ell'$, but the string of length $\ell'$ 
in $(\nu,J)^{(k)}$ is singular which implies that $\ell^{(k)}\le \ell'<\ell$. Hence 
there must exist $\ell'<i<\ell$ such that $m_i^{(k-1)}(\nu)>0$ and $\ell^{(k-1)}=i$.
By \eqref{eq:dp ineq} the same five cases as in Remark~\ref{remark:new tab}
occur as possibilities for the letters $i$ and $i+1$ in columns $k$ and $k+1$ of $t$. 
In cases~\ref{l:3}, \ref{l:4} and case~\ref{l:5} if $m_i^{(k-1)}(\nu)=2$,
we have $m_i^{(k+1)}(\nu)=0$. Replace $i$ in column $k+1$ by $i+1$ in $t$ to get
a new tableau $t'$. In all other cases $m_i^{(k-1)}(\nu)=1$; replace the letter $i+1$
in column $k$ by $i$ to obtain $t'$.
The replacement $t\mapsto t'$ yields $\Dp_i^{(k)}(t')>0$ in all cases.
The change of lower bound $M_i^{(k-1)}(t')=M_i^{(k-1)}(t)+1$ in cases \ref{l:1},
\ref{l:2} and \ref{l:5} when $m_i^{(k-1)}(\nu)\neq 2$ will not cause any
problems since $m_i^{(k-1)}(\nu)=1$ so that after the application
of $\delta$ there is no part of length $i$ in the $(k-1)$-th rigged partition. 
Then again repeated application of Remark~\ref{remark:new tab} achieves
$\Dp_{\ell-1}^{(k)}(t'')>0$ for some tableau $t''$, so that the problematic
case does not occur. 

Let $t''$ be the tableau we constructed so far. Note that in all constructions 
above,  either a letter $i+1$ in column $k$ is changed to $i$, or a letter $i$ 
in column $k+1$ is changed to $i+1$. In the latter case $i+1\le \ell\le 
|\nu^{(k)}|\le c_k$. Hence $t''$ satisfies the constraint that 
$t''_{i,k}\in\{1,2,\ldots,c_{k-1}\}$ for all $i,k$.

Now let  $\tb=\D_r(t'')$. We know $\tb \in \A(\lab')$.
We will show that  the parts of $\Jb$ lie between the corresponding lower 
bound with respect to $\tb \in \A(\lab')$ and the vacancy number. 

If $t''_{1,r}<c_{r-1}$ then by Lemma~\ref{lem:biggest part} $M_i^{(k)}(\tb)\le
M_i^{(k)}(t'')$  for all $k$ and $i$ such that $m_i^{(k)}(\nub)>0$.  Hence 
by Lemma ~\ref{lem:convex}  we have that  $(\nub,\Jb)$ is admissible with 
respect to  $\tb$. 

Let  $t''_{1,r}=c_{r-1}$. Then there exists $j$ as in the definition of $\D_r$. 
We claim that 
\begin{enumerate}
\item[(i)] $m_i^{(r-1)}(\nu)=0$ for $i>c_{r-1}-j$ and $m_{c_{r-1}-j}^{(r-1)}(\nu)\le 1$.
\item[(ii)] If $m_{c_{r-1}-j}^{(r-1)}(\nu)=1$, then $\ell^{(r-1)}=c_{r-1}-j$. 
\end{enumerate}
Note that $M_i^{(r-1)}(\tb)=M_i^{(r-1)}(t'')+1$ for $c_{r-1}-j\le i <c_{r-1}$ and 
$M_i^{(k)}(\tb)\le M_i^{(k)}(t'')$ for all other $k$ and $i$ such that 
$m_i^{(k)}(\nub)>0$. Hence if  the claim is true  using  Lemma~\ref{lem:biggest part}  
we  have $M_i^{(k)}(\tb)\le M_i^{(k)}(t'')$  for all $k$ and $i$ such that 
$m_i^{(k)}(\nub)>0$. Therefore by Lemma~\ref{lem:convex} we have that  $(\nub,\Jb)$ 
is admissible with respect to $\tb$.

It remains to prove the claim. Note that if $|\nu^{(r-1)}|<c_{r-1}-j$ then our claim 
is trivially true. Let $|\nu^{(r-1)}|=c_{r-1}-k$ for some $0\le k \le j$. If all parts of
$\nu^{(r-1)}$ are strictly less than $c_{r-1}-j$, again our claim is trivially true. 
Let the largest part in $\nu^{(r-1)}$ be $c_{r-1}-p\ge c_{r-1}-j$ for some $k\le p\le j$. 
Let $a$ be the largest part in $\nu^{(r)}$. 

First suppose $a> c_{r-1}-p$ and $a=c_r-q$ for some
$0\le q <c_r$. Then $a=c_{r-1}-(\la_r+q)$ which implies that
\begin{equation*}
M_a^{(r)}(t'')\ge-(c_r-\la_r-q)+(c_{r+1}-q)=  \la_r-\la_{r+1}.
\end{equation*}
This means $p_a^{(r)}(\nu)\le M_a^{(r)}(t'')$ since $p_b^{(r)}(\nu)\ge p_a^{(r)}(\nu)$ 
for all $b\ge a$ and $p_b^{(r)}=\la_r-\la_{r+1}$ for large $b$. If 
$p_a^{(r)}(\nu)< M_a^{(r)}(t'')$,  it contradicts that
$p_a^{(r)}(\nu)\ge M_a^{(r)}(t'')$. If $p_a^{(r)}(\nu)= M_a^{(r)}(t'')$,  it contradicts the 
fact that $r=\rk(\nu,J)$ since we get a singular part of length $a$ in $\nu^{(r)}$
which is larger than the largest part in $\nu^{(r-1)}$. 
Therefore $a> c_{r-1}-p$ is not possible. 

Hence $a\le c_{r-1}-p$. Using Remark~\ref{remark:rigsize} we get,
\begin{equation}\label{eq:max vacancy}
\begin{split} 
p_a^{(r)}(\nu)&=Q_a(\nu^{(r-1)})-2|\nu^{(r)}|+Q_a(\nu^{(r+1)}) +\sum_{i\ge 1}
\min(a,i)L_i^{(r)}\\
&\le  a+p-k-2|\nu^{(r)}|+|\nu^{(r+1)}| +\sum_{i\ge 1}\min(a,i)L_i^{(r)}\\
&=a+p-2\la_{r+1}-c_{r+1}-\sum_{i\ge 1}\max(s_i-a,0)\delta_{r_i,r}.
\end{split}
\end{equation}    
Since $ p_a^{(r)}(\nu)\ge M_a^{(r)}(t'')\ge -\la_{r+1}$ we get 
\begin{equation*}
c_{r}-(p-\sum_{i\ge 1}\max(s_i-a,0)\delta_{r_i,r})\le a\le c_r.
\end{equation*}
Hence $a=c_r-q$ for $0\le q\le p- \sum_{i\ge 1}\max(s_i-a,0)\delta_{r_i,r}$.
Then from \eqref{eq:max vacancy} with $a=c_r-q$ we get
\begin{equation}\label{eq:ub vacancy}
 p_a^{(r)}(\nu)\le p-q-\la_{r+1}- \sum_{i\ge 1}\max(s_i-a,0)\delta_{r_i,r}
\le \la_r-\la_{r+1},
\end{equation}
where we used that $0\le p-q \le \la_r$ which follows from $a=c_r-q\le c_{r-1}-p$.

If $a>c_{r-1}-j$, as in the case $a>c_{r-1}-p$ we have
$$ M_a^{(r)}(t'')\ge -(c_r-\la_r-q)+(c_{r+1}-q)= \la_r-\la_{r+1}\ge p_a^{(r)}(\nu).$$ 
Hence we get a contradiction unless $p_a^{(r)}(\nu)= M_a^{(r)}(t'')$. By 
\eqref{eq:ub vacancy} and the fact that $0\le p-q\le \la_r$ we know  
$p_a^{(r)}(\nu)=\la_r-\la_{r+1}$ happens only when $p-q=\la_r$ and 
$\sum_{i\ge 1}\max(s_i-a,0)\delta_{r_i,r}=0$. This means the largest
part in $\nu^{(r-1)}$ is of length $c_{r-1}-p=c_r-q=a$. Since we have a 
singular string of length $a$ in $\nu^{(r)}$  this contradicts the fact that $r=\rk(\nu,J)$. 

If $a\le c_{r-1}-j$ then $ M_a^{(r)}(t'')\ge -(c_r-j)+(c_{r+1}-q)= j-q-\la_{r+1}\ge 
p_a^{(r)}(\nu)$ 
because of \eqref{eq:ub vacancy} and the fact that $j\ge p$. Again we get a contradiction 
unless $p_a^{(r)}(\nu)= M_a^{(r)}(t'')$. But this happens only when 
$p_a^{(r)}(\nu)=j-q-\la_{r+1}$ which gives $p=j$ because $p_a^{(r)}(\nu)$ attains 
the right hand side of \eqref{eq:ub vacancy}.  This means the largest part in 
$\nu^{(r-1)}$ is $c_{r-1}-j$. Furthermore, for large $i$ we have
$p_i^{(r)}=\la_r-\la_{r+1}\ge j-q-\la_{r+1}+(c_{r-1}-j-a)=\la_r-\la_{r+1}$ which
shows that besides $c_{r-1}-j$ all parts in $\nu^{(r-1)}$ have to be less 
than or equal to $a$. But the part of length $a$ in $\nu^{(r)}$ is singular, so 
we have to have $c_{r-1}-j>a$ and $\ell^{(r-1)}=c_{r-1}-j$ else it will 
contradict the fact that $r=\rk(\nu,J)$. This proves our claim.   

Hence $(\nub,\Jb)$ is admissible with respect to $\tb \in \A(\lab')$ and therefore
$\delta$ is well-defined.           
\end{proof}

\begin{example}
Let $L$ be the multiplicity array of $B=(B^{1,1})^{\otimes 4}$ and $\la=
(0,1,0,1,2)$. Let 
\begin{equation*}
(\nu,J)= \yngrc(3,-1,1,2) \quad \yngrc(2,0,1,0) \quad \yngrc(2,-1,1,-1) \quad
\yngrc(2,-1)\in \RC(L,\la).
\end{equation*} 
Let $t=\young(4433,3222,211,1)$ be the corresponding lower bound 
tableau. Then 
\begin{equation*}
\delta(\nu,J)= \yngrc(3,-1) \quad \yngrc(2,0) \quad \yngrc(2,-1) \quad
\yngrc(1,-1).
\end{equation*} 
Note that in this example $\ell=\ell^{(4)}=2$ and it satisfies \eqref{eq:problem}
with $k=4$. Also $\Dp_{\ell-2}^{(4)}(t)+\Dp_{\ell}^{(4)}(t)=0$ with 
$\Dp_{i}^{(4)}(t)=0$ for all $0\le i\le \ell$. Since $m_1^{(3)}(\nu)=1$ and
$2\in t_{.,4}$ this is an example where we get the new tableau $t'$ by replacing
the $2\in t_{.,4}$ by $1$ and then the corresponding lower 
bound tableau for $\delta(\nu,J)$ is  $\D_5(t')=\young(3221,211,1)$.
\end{example}

\begin{proof}[Proof of Proposition~\ref{prop:inv delta}]
Similar to Proposition~\ref{prop:delta} we need to show that for $(\nub,\Jb)\in
\RC(\overline{L},\lab)$ we have $\delta^{-1}(\nub,\Jb)=(\nu,J)\in \RC(L,\la)$ 
where $\la=\lab+\epsilon_r$. Clearly $\la$
has nonnegative parts, so it suffices to show that $(\nu, J)$ is admissible which 
means that the parts of $J$ lie between the corresponding lower bound with 
respect to some $t \in \A (\la')$ and the vacancy number. Let 
$\tb \in \A(\lab')$ be a tableau such that $(\nub, \Jb)$ is admissible with 
respect to $\tb$.  By similar argument as in the proof of Propostion~\ref{prop:delta} 
the only problematic case occurs when
\begin{equation}\label{eq:inv problem}
m_{s+1}^{(k)}(\nub)=0, \quad \Dp_{s+1}^{(k)}(\tb)=0, \quad s<s^{(k+1)}
\quad \text{and $ s$ finite}
\end{equation}
where $s=s^{(k)}$.

Assume that $\Dp_{s}^{(k)}(\tb)+\Dp_{s+2}^{(k)}(\tb)\ge 1$ and \eqref{eq:inv problem}
holds. By Remark~\ref{remark:new tab} with $i=s+1$  there exists a new tableau
$\tb'$ such that $\Dp_{s+1}^{(k)}(\tb')>0$ so that the problematic case is avoided.

Hence assume that $\Dp_{s}^{(k)}(\tb)+\Dp_{s+2}^{(k)}(\tb)=0$ and \eqref{eq:inv problem}
holds. Let $s'>s$ be minimal such that $m_{s'}^{(k)}(\nub)>0$. If no such
$s'$ exists, set $s'=\infty$. 

Suppose that there exists $s'>j>s$ such that $\Dp_{j+1}^{(k)}(\tb)>0$.
Let $i$ be the minimal such $j$. Then by Remark~\ref{remark:new tab}
we can find a new tableau $\tb'$ such that $\Dp_i^{(k)}(\tb')>0$ and 
$(\nub,\Jb)$ is admissible with respect to $\tb'$.
Repeating the argument we can achieve $\Dp_{s+1}^{(k)}(\tb'')>0$ for
some new tableau $\tb''$, so that the problematic case does not occur.

Hence we are left to consider the case $\Dp_i^{(k)}(\tb)=0$ for all
$s'\ge i \ge s$.  First let us suppose $k<r-1$. If $m_i^{(k+1)}(\nub)=0$ for all $s'>i>s$, 
then by the similar arguments as in the proof of Proposition ~\ref{prop:delta} we arrive 
at a contradiction since $s^{(k+1)}\ge s'$, but the string of length $s'$ 
in $(\nub,\Jb)^{(k)}$ is singular which implies that $s^{(k+1)}>s^{(k)}\ge s'>s$. Hence 
there must exist $s'>i>s$ such that $m_i^{(k+1)}(\nub)>0$ and $s^{(k+1)}=i$.
By \eqref{eq:dp ineq} the same five cases as in Remark~\ref{remark:new tab}
occur as possibilities for the letters $i$ and $i+1$ in columns $k$ and $k+1$ of $\tb$. 
In cases~\ref{l:1}, \ref{l:2} and case~\ref{l:5} if $m_i^{(k+1)}(\nub)=2$,
we have $m_i^{(k-1)}(\nub)=0$. Replace $i+1$ in column $k$ by $i$ in $\tb$ to get
a new tableau $\tb'$. In all other cases $m_i^{(k-1)}(\nub)=1$; replace the letter $i$
in column $k+1$ by $i+1$ to obtain $\tb'$.
The replacement $\tb\mapsto \tb'$ yields $\Dp_i^{(k)}(\tb')>0$ in all cases.
The change of lower bound $M_i^{(k+1)}(\tb')=M_i^{(k+1)}(\tb)+1$ in cases \ref{l:3},
\ref{l:4} and \ref{l:5} when $m_i^{(k+1)}\neq 2$ will not cause any
problems since $m_i^{(k+1)}=1$ so that after the application
of $\delta^{-1}$ there is no part of length $i$ in the $(k+1)$-th rigged partition. 
Then again repeated application of Remark~\ref{remark:new tab} achieves
$\Dp_{s+1}^{(k)}(\tb'')>0$ for some tableau $\tb''$, so that the problematic
case does not occur. 

Now let us consider the case $k=r-1$. Note that $s'=\infty$ here. Else $s^{(r-1)}>s$, 
a contradiction. So, $\Dp_i^{(r-1)}(\tb)=0$ for $i>s$ which implies 
$m_i^{(r-1)}(\nub)=0$ for $i>s$, else $s^{(r-1)}>s$.  Then
by~\eqref{eq:dp ineq} with $i\ge s+1$ and $k=r-1$ we have 
\begin{equation} \label{eq:S ineq}
\begin{split}
-\chi(i\in \tb_{\cdot,r-1})+&\chi(i\in \tb_{\cdot,r})+\chi(i+1\in \tb_{\cdot,r-1})
-\chi(i+1\in \tb_{\cdot,r})\\
&\ge m_{i}^{(r-2)}(\nub)+m_{i}^{(r)}(\nub)\ge 0.
\end{split}
\end{equation}

If  $s+1\in \tb_{.,r}$  by \eqref{eq:S ineq} with $i=s+1$ there are seven choices for the 
letters $s+1$ and $s+2$ in columns $r-1$ and $r$ of $\tb$.
\begin{enumerate}
\item \label{s:1} $s+1$ in both columns $r-1$ and $r$;
\item \label{s:2} Both $s+1,s+2$ in  column $r$; 
\item \label{s:3} Both $s+1,s+2$ in columns $r-1,r$;
\item \label{s:4}  $s+1$ in  columns $r-1,r$ and $s+2$ in column $r-1$; 
\item \label{s:5} $s+1$ in column  $r$;
\item \label{s:6} $s+1$ in  column $r$ and $s+2$ in columns $r-1,r$;
\item \label{s:7} $s+1$ in  column $r$ and $s+2$ in column $r-1$.
\end{enumerate} 
First note that by \eqref{eq:S ineq} $m_{s+1}^{(r-2)}(\nub)=m_{s+1}^{(r)}(\nub)=0$ 
for cases \ref{s:1}, \ref{s:2} and \ref{s:3}. For case \ref{s:4} we have 
$m_{s+1}^{(r)}(\nub)=0$ again, else  $p_{s+1}^{(r-1)}(\tb)>p_{s}^{(r-1)}(\tb)=
M_{s}^{(r-1)}(\tb)=M_{s+1}^{(r-1)}(\tb)$,
contradiction to $\Dp_{s+1}^{(r-1)}(\tb)=0$. In cases~\ref{s:5} and~\ref{s:6}  
either $m_{s+1}^{(r)}(\nub)=0$ or $m_{s+1}^{(r-2)}(\nub)=0$ by \eqref{eq:S ineq}. 
When  $m_{s+1}^{(r-2)}(\nub)=0$ and  $m_{s+1}^{(r)}(\nub)>0$ in case~\ref{s:5} 
we have $m_i^{(r-2)}(\nub)=0$ for all $i>s+1$, else 
$p_{s+1}^{(r-1)}(\nub)\ge p_{s}^{(r-1)}(\nub)+2=
M_{s}^{(r-1)}(\tb)+2\ge M_{s+1}^{(r-1)}(\tb)-1+2>M_{s+1}^{(r-1)}(\tb)$, 
a contradiction. In case~\ref{s:7} by the same string of inequalities 
either $m_{s+1}^{(r)}(\nub)=0$ or $m_{s+1}^{(r-2)}(\nub)=0$. 

When $m_{s+1}^{(r)}(\nub)=0$ we construct a new tableau $\tb'$ from $\tb$ 
by replacing $s+1$ in column $r$ by the smallest number $i>s+1$ that  
does not appear in column $r$ of $\tb$. The effect of this change is 
$M_{s+1}^{(r)}(\tb') =M_{s+1}^{(r)}(\tb)+1$ and $M_{s+1}^{(r-1)}(\tb') =
M_{s+1}^{(r-1)}(\tb)-1$. Since $m_{s+1}^{(r)}(\nub)=0$ the first change
does not create any problem.  When $m_{s+1}^{(r)}(\nub)>0$  in cases~\ref{s:6} 
and~\ref{s:7} we change the $s+2$ in column $r-1$ to 
$s+1$. The effect of this replacement is $M_{s+1}^{(r-2)}(\tb') 
=M_{s+1}^{(r-2)}(\tb)+1$ and $M_{s+1}^{(r-1)}(\tb') =M_{s+1}^{(r-1)}(\tb)-1$. 
Since $m_{s+1}^{(r-2)}(\nub)=0$ there is no problem. When 
$m_{s+1}^{(r)}(\nub)>0$ in case~\ref{s:5} we replace the smallest 
$\tb_{j,r-1}>s+1$ by $s+1$. This has the effect that 
$M_{i}^{(r-2)}(\tb') =M_{i}^{(r-2)}(\tb)+1$ for $s+1\le i< \tb_{j,r-1}$. Since 
we have $m_i^{(r-2)}=0$ for all $i\ge s+1$ we do not have any problem. 
In all cases, replacing $\tb$ by $\tb'$ the problematic case~\eqref{eq:inv problem} 
is avoided and we have $\Dp_i^{(k)}(\tb')\ge 0$ for 
all other $i,k$ such that $m_i^{(k)}(\nub)>0$. 

Let us consider the case $s+1 \not \in \tb_{.,r}$. Note that
$M_{s}^{(r-1)}(\tb)\ge M_{s+1}^{(r-1)}(\tb)$. 
We have $m_i^{(r)}(\nub)=0=m_i^{(r-2)}(\nub)$ for all 
$i>s$, else $p_{s+1}^{(r-1)}(\nub)>p_{s}^{(r-1)}(\nub)=M_{s}^{(r-1)}(\tb)\ge 
M_{s+1}^{(r-1)}(\tb)$, contradiction to $\Dp_{s+1}^{(r-1)}(\tb)=0$. Using   
\eqref{eq:S ineq}   for $i=s+1, k=r-1$ we have four possible cases for the 
choice of the letters $s+1$ and $s+2$ in columns $r-1$ and $r$ of $\tb$. 
 \begin{enumerate}
 \item \label{i:1} $s+2$ in column $r-1$;
 \item \label{i:2} $s+2$ in columns $r-1$ and $r$;
 \item \label{i:3} $s+1$ and $s+2$ in column $r-1$;
 \item \label{i:4} no $s+1,s+2$ in both columns $r-1$ and $r$.
\end{enumerate}
We first argue that case~\ref{i:3} cannot occur.
Suppose case~\ref{i:3} holds. Then $M_{s+1}^{(r-1)}(\tb)= M_{s}^{(r-1)}(\tb)-1$ and
$M_{s+2}^{(r-1)}(\tb)= M_{s+1}^{(r-1)}(\tb)-1$. But we also have 
$\Dp_i^{(r-1)}(\tb)=0$ for $i>s$ and $m_{i}^{(r-1)}(\nub)= m_{i}^{(r-2)}(\nub)= 
m_{i}^{(r)}(\nub)$ for $i>s$. Note that $\Dp_i^{(r-1)}(\tb)=0$ implies that 
$p_{s+2}^{(r-1)}(\nub)= p_{s+1}^{(r-1)}(\nub)-1=p_s^{(r-1)}(\nub)-2$. On the other hand
$m_{i}^{(r-1)}(\nub)= m_{i}^{(r-2)}(\nub)= m_{i}^{(r)}(\nub)$ implies that 
$p_{s+2}^{(r-1)}(\nub)\ge p_s^{(r-1)}(\nub)$ and $p_{s+1}^{(r-1)}(\nub)\ge 
p_s^{(r-1)}(\nub)$ which yields a contradiction.

In cases~\ref{i:1} and~\ref{i:2} we replace the letter $s+2$ in column $r-1$ to 
$s+1$ to get a new tableau $\tb'$. The change from $\tb$ to $\tb'$ yields 
$\Dp_{s+1}^{(r-1)}(\tb')>0$ without any other change. 
In case~\ref{i:4} if there exists $\tb_{j,r-1}>s+2$ for some $j$ then we 
replace the smallest such $\tb_{j,r-1}$ by $s+1$ to construct $\tb'$. Then 
again we get $\Dp_{s+1}^{(r-1)}(\tb')>0$ without any other change since 
$m_i^{(r-2)}(\nub)=0$ for all $i>s$. On the other hand if $\tb_{1,r-1}\le s$ then  
$\cb_{r-1} \le s\le |\nub^{(r-1)}|\le \cb_{r-1}$ implies $\tb_{1,r-1}=s$. 
Note that $\tb_{1,r-2}\ge s$. Here  we will avoid the problematic 
case \eqref{eq:inv problem} by constructing a new 
tableau $t\in \A(\la')$. Let 
\begin{equation}\label{eq:sp t}
t_{i,k} = \begin{cases}
 \cb_1+1 & \text{for $k=1=i$}\\
 \cb_{k-1} +1& \text{for $2\le k\le r-2$ and $i=1$},\\
 s+1 & \text{for  $k=r-1$ and $i=1$},\\ 
 \tb_{i-1,k} & \text{for $1\le k\le r-1$ and $1< i\le \cb_k$},\\
 \tb_{i,k} & \text{for $r\le k\le n$ and $1\le i \le \cb_k$.}
\end{cases}
\end{equation} 
Note that $c_k=\cb_k+1$ for $1\le k \le r-1$
and $c_k=\cb_k$ for $r\le k \le n$. Clearly $t_{i,k}\in \{1,2,\ldots,c_{k-1}\}$ 
for all $i,k$. Column-strictness of $t$ follows since $\tb_{1,1}< \cb_1+1$ 
and $\tb_{1,k}<\cb_k+1 \le \cb_{k-1}+1$ for
$2\le k\le r-1$ and $s+1>\tb_{1,r}$ .  Hence $t \in \A(\la')$. 
Note that we have $M_{s+1}^{(r-1)}(t)=  M_{s+1}^{(r-1)}(\tb)-1<  
p_{s+1}^{(r-1)}(\nub)$, so the problematic case~\eqref{eq:inv problem} is
avoided. The fact that $(\nu,J)$ is admissible with respect to $t$ is shown
later. 

Let us now define $t\in \A(\la')$ in all other cases. 
Let $\tb''\in \A(\lab')$ be the tableau we constructed from $\tb$ so far 
except in the last case. Note that in all constructions above,  either a 
letter $i+1$ in column $k$ is changed to $i$, or 
a letter $i$ in column $k+1$ is changed to $i+1$. In the latter case 
$m_i^{(k+1)}=0$ means $i+1\le s^{(k+1)}\le |\nu^{(k+1)}|\le \cb_{k+1}\le 
\cb_k$. Hence $\tb''$ satisfies the constraint that 
$\tb''_{i,k}\in\{1,2,\ldots,\cb_{k-1}\}$ for all $i,k$. 

Let us define a new tableau $t$ from $\tb''$ in the following way:
\begin{equation}\label{eq:t}
t_{i,k} = \begin{cases}
 \cb_1+1 & \text{for $k=1=i$}\\
 \cb_{k-1} +1& \text{for $2\le k\le r-1$ and $i=1$},\\ 
 \tb''_{i-1,k} & \text{for $1\le k\le r-1$ and $1< i\le \cb_k$,}\\
 \tb''_{i,k} & \text{for $r\le k\le n$ and $1\le i \le \cb_k$.}
\end{cases}
\end{equation}    
Similarly as in \eqref{eq:sp t} we have  $t \in \A(\la')$. 

Next we show that $(\nu,J)$ is admissible with respect to $t$, that is,
the parts of $J$ lie between the corresponding lower bound with 
respect to $t \in \A(\la')$ and the vacancy number.  
Note that $s^{(k)}+1 \le |\nu^{(k)}| \le c_k \le c_{k-1}$. We distinguish the
three cases $s^{(k)}+1<c_k$, $s^{(k)}+1=c_k=c_{k-1}$ and $s^{(k)}+1=c_k<c_{k-1}$.

If $s^{(k)}+1<c_k$ for all $1\le k \le r-1$, then $M_i^{(k)}(t)=M_i^{(k)}(\tb'')$ for 
all $i,k$ such that $m_i^{(k)}(\nu)>0$.
If $s^{(k)}+1=c_{k-1}$ for some $1\le k \le r-2$, then $M_{s^{(k)}+1}^{(k)}(t)=
M_{s^{(k)}+1}^{(k)}(\tb'')$ since $c_{k-1}\ge c_{k}$. Also if
$s^{(r-1)}+1=c_{r-2}$, then $M_{s^{(r-1)}+1}^{(r-1)}(t)=
M_{s^{(r-1)}+1}^{(r-1)}(\tb'')-1$. In both cases $(\nu,J)$ is admissible 
since $M_i^{(k)}(t)\le M_i^{(k)}(\tb'')$ for all $i,k$ such that $m_i^{(k)}(\nu)>0$.

Now suppose $s^{(k)}+1=c_{k} < c_{k-1}$  for some $1\le k < r-1$. Then 
$M_{s^{(k)}+1}^{(k)}(t)=M_{s^{(k)}+1}^{(k)}(\tb'')+1$. Suppose $k$ is 
minimal satisfying this condition. Note that in this situation, $s^{(k)}= 
c_k-1=\cb_k$. This means $|\nub^{(k)}|=\cb_k$ which implies by 
definition of $|\nub^{(k)}|$ that $|\nub^{(a)}|=\cb_a$ for $a\ge k$.
Using this we get 
$$ \cb_k=s^{(k)}\le s^{(k+1)}\le\cdots \le s^{(a)}\le \cdots \le s^{(r-1)}\le |
\nub^{(r-1)}|=\cb_{r-1}\le \cb_k.$$ 
This implies  $\cb_a=s^{(a)}=s^{(a+1)}=\cb_{a+1}$ for all $k\le a \le r-2$. 
When  $s^{(a)}=s^{(a+1)}$ we have   $p_{s^{(a)}+1}^{(a)}(\nu)=
p_{s^{(a)}+1}^{(a)}(\nub)$. Hence we only need to worry  
when $\Dp_{s^{(k)}+1}^{(k)}(\tb'')= 0$.  Let $\ell$ be the largest  
part in $\nub^{(k-1)}$. If $\ell>s^{(k)}$ then by definition 
$p_{s^{(k)}+1}^{(k)}(\nub)>p_{s^{(k)}}^{(k)}(\nub)$. 
But we have $M_{s^{(k)}}^{(k)}(\tb'')\ge M_{s^{(k)}+1}^{(k)}(\tb'')$, hence  
$\Dp_{s^{(k)}+1}^{(k)}(\tb'')> 0$. Suppose $\ell\le s^{(k)}$, then  
$p_{s^{(k)}+1}^{(k)}(\nub)\ge p_{s^{(k)}}^{(k)}(\nub)$ since 
$m_i^{(k)}(\nub) =0$ for $i>s^{(k)}$. If $s^{(k)}+1\in \tb''_{.,k}$ 
then $M_{s^{(k)}}^{(k)}(\tb'')= M_{s^{(k)}+1}^{(k)}(\tb'')+1$ and we 
get $\Dp_{s^{(k)}+1}^{(k)}(\tb'')> 0$.  If  $s^{(k)}+1\not \in \tb''_{.,k}$
then there exists $\tb''_{j,k}>s^{(k)}+1$ for some $j$ and we replace 
the  smallest such $\tb''_{j,k}$  by $s^{(k)}+1$ to get a new tableau $t'$ 
from $t\in \A(\la')$ . This has the effect that  $M_{s^{(k)}+1}^{(k)}(t')=
M_{s^{(k)}+1}^{(k)}(t)-1=M_{s^{(k)}+1}^{(k)}(\tb'')$ so that 
$\Dp_{s^{(k)}+1}^{(k)}(t')\ge 0$. 

This proves  that  $(\nu,J)$ is admissible with respect to 
$t$ or $t' \in\A(\la') $. Hence $\delta^{-1}$ is well-defined.
\end{proof}

\begin{example}
Let $\overline{L}$ be the multiplicity array of $B=(B^{1,1})^{\otimes 4}$ and $\lab=
(0,1,1,1,1)$. Let 
\begin{equation*}
(\nub,\Jb)= \yngrc(3,-1,1,1) \quad \yngrc(2,-1,1,0) \quad \yngrc(1,-1,1,-1) \quad
\yngrc(1,0)\in \RC(\overline{L},\lab).
\end{equation*} 
Let $\tb=\young(4432,321,21,1)$ be the corresponding lower bound 
tableau. Then with $r=3$,
\begin{equation*}
\delta^{-1}(\nub,\Jb)= \yngrc(3,-1,1,1,1,1) \quad \yngrc(3,-1,1,1) \quad \yngrc(1,-1,1,-1)
 \quad \yngrc(1,0).
\end{equation*} 
Note that in this example we have $k=r-1=2$ and $s=s^{(2)}=2$ which satisfies 
\eqref{eq:inv problem}. Also $s+1=3\in \tb_{.,r}$, hence this is the situation when $k=r-1$ 
in \eqref{eq:inv problem} with $\Dp_i^{(r-1}(\tb)=0$ for all $i>s$ and since 
$s+1\in \tb_{.,r}$ this is case~\ref{s:7} discussed in the proof. So we get the 
corresponding 
lower bound tableau for $(\nu,J)$ by replacing $3\in \tb_{.,r}$ by $4$ and then doing the 
construction defined in \eqref{eq:t}. The lower bound tableau we get is 
$\young(5542,441,32,21,1)$.
\end{example}

\section{Proof of Proposition~\ref{prop:bij}}\label{appn:phi}
In this section a proof of Proposition~\ref{prop:bij} is given stating that
the map $\Phi$ of Definition~\ref{def:bij} is a well-defined bijection.

The proof proceeds by induction on $B$ using the fact that 
it is possible to go from $B=B^{r_k,s_k}\otimes B^{r_{k-1},s_{k-1}}
\otimes \cdots \otimes B^{r_1,s_1}$ to the 
empty crystal via successive application of $\lh$, $\ls$ and $\lb$. 
Suppose that $B$ is the empty crystal. Then both sets $\Path(B,\la)$ and 
$\RC(L,\la)$ are empty unless $\la$ is the empty partition, in which case 
$\Path(B,\la)$ consists of the empty partition and $\RC(L,\la)$ consists of 
the empty rigged configuration. In this case $\Phi$ is the unique
bijection mapping the empty partition to the empty rigged configuration.

Consider the commutative diagram~\eqref{bij:1} of Definition~\ref{def:bij}.
By induction
\begin{equation*}
 \Phi: \displaystyle{\bigcup_{\mu\in\lm} \Path(\lh(B),\mu)} \longrightarrow  
\displaystyle{\bigcup_{\mu\in\lm} \RC(\lh(L),\mu)}
\end{equation*}
is a bijection. By Propositions~\ref{prop:delta} and~\ref{prop:inv delta}
$\delta$ is a bijection, and by definition it is clear that $\lh$ is a bijection
as well. Hence $\Phi=\delta^{-1} \circ \Phi \circ \lh$ is a well-defined bijection.

Suppose that $B=B^{r,1} \otimes B'$ with $r\ge2$. By induction $\Phi$ is a bijection
for $\lb(B)=B^{1,1}\otimes B^{r-1,1} \otimes B'$. Hence to prove 
that~\eqref {bij:3} uniquely determines $\Phi$ for $B$ it suffices to show that 
$\Phi$ restricts to a bijection between the image of 
$\lb:  \Path(B,\la)\longrightarrow \Path(\lb(B),\la)$ and the
image of $\rclb : \RC(L,\la) \longrightarrow \RC(\lb(L),\la)$.
Let $b=\begin{array}{|c|} \hline b_r\\ \hline \end{array} \otimes 
\begin{array}{|c|} \hline b_1\\ \hline \vdots \\ \hline b_{r-1}\\ \hline \end{array} 
\otimes b' \in \Path(\lb(B),\la)$ with $b_{r-1}<b_r$. Let $(\nu,J)=\Phi(b)$ which
is in $\RC(\lb(L),\la)$. We will show that $(\nu,J)^{(a)}$ has a 
singular string of length one for $1\le a\le r-1$.  

By induction we know  for $(\nub,\Jb)=\Phi(\overline{b})$ where $\overline{b}
=\begin{array}{|c|} \hline b_{r-1}\\ \hline \end{array} \otimes 
\begin{array}{|c|} \hline b_1\\ \hline \vdots \\ \hline b_{r-2}\\ \hline \end{array} 
\otimes b' \in \lb(B^{r-1,1}\otimes B')$ with $b_{r-2}<b_{r-1}$, $(\nub,\Jb)^{(a)}$ 
has a singular string of length one for $1\le a\le r-2$.  Let $\overline{b}'= 
\begin{array}{|c|} \hline b_1\\ \hline \vdots \\ \hline b_{r-1}\\ \hline \end{array} 
\otimes b'$ and $(\nub',\Jb')=\Phi(\overline{b}')$. This "unsplitting" on the
rigged configuration  side removes the singular string of length one from 
$(\nub,\Jb)^{(a)}$  for $1\le a \le r-2$ yielding $(\nub',\Jb')$.

Let  $\s^{(a)}$ be the length of the selected strings by $\delta^{-1}$ 
associated with $b_{r-1}$. Note that $\s^{(a)}=0$  for $1\le a \le r-2$. 
Now let $s^{(a)}$ be the selected strings by $\delta^{-1}$ 
associated with $b_r$. Since $b_{r-1}<b_r$ we have by 
construction that $s^{(a+1)}\le \s^{(a)}$. In
particular $s^{(r-1)}\le \s^{(r-2)}=0$ and therefore,
$s^{(r-1)}=0$. This implies that $s^{(a)}=0$ for $1\le a\le r-1$. 
Hence $(\nu,J)^{(a)}$ has a singular 
string of length one for $1\le a \le r-1$. 

Conversely, let $(\nu,J)\in \rclb(\RC(L,\la))$, that is, $(\nu,J)^{(a)}$ has 
singular string of length one for $1\le a \le r-1$. Let $b=\Phi^{-1}(\nu,J)=
\begin{array}{|c|} \hline b_r\\ \hline \end{array} \otimes 
\begin{array}{|c|} \hline b_1\\ \hline \vdots \\ \hline b_{r-1}\\ \hline \end{array} 
\otimes b' \in \Path(\lb(B),\la)$. We want to show that $b_{r-1}<b_r$. 
Let $(\nub,\Jb)=\delta(\nu,J)$ and $\ell^{(a)}$ be the length of the 
selected string in $(\nu,J)^{(a)}$ by $\delta$. Then $\ell^{(a)}=1$
for $1\le a \le r-1$ and the change of vacancy numbers from $(\nu,J)$ 
to $(\nub,\Jb)$ is given by
\begin{equation}\label{eq:l change in vac}
p_i^{(a)}(\nub)= p_i^{(a)}(\nu)-\chi(\ell^{(a-1)}\le i < \ell^{(a)})
 +\chi( \ell^{(a)}\le i < \ell^{(a+1)}).
\end{equation}
This implies that $(\nub,\Jb)^{(r-1)}$ has no singular string of 
length less than $\ell^{(r)}$ since $\ell^{(r-1)}=1$. Let 
$(\nub',\Jb')=\rclb(\nub,\Jb)$. Denote by $\ellb^{(a)}$ the length of the singular 
string selected by $\delta$ in $(\nub',\Jb')^{(a)}$.  
Then by induction $\ellb^{(a)}=1$ for $1\le a \le r-2$ and by 
\eqref{eq:l change in vac} we get $\ellb^{(a)}\ge \ell^{(a+1)}$ for 
$a\ge r-1$. Therefore $\ellb^{(a)}\ge \ell^{(a+1)}$ for all $1\le a\le n$. 
Hence $b_{r-1}<b_r$. This proves that $\Phi$ in~\eqref{bij:3} is
uniquely determined.

Let us now consider the case $B=B^{r,s}\otimes B'$ where $s\ge 2$. 
Any map $\Phi$ satisfying \eqref{bij:2} is injective by definition and 
unique by induction. To prove the existence and surjectivity it 
suffices to prove that  bijection $\Phi$ maps the image of $\ls: \Path(B,\la)
\longrightarrow \Path(\ls(B),\la)$ to the image of $\rcls: \RC(L,\la) 
\longrightarrow \RC(\ls(L),\la)$. 
Let $b=c_1\otimes c \otimes b' \in \ls(\Path(B,\la))$ 
where $c=c_2c_3\cdots c_{s}$ and $c_i$
denotes the $(i-1)$-th column of $c\in B^{r,s-1}$. Let $c_1= \begin{array}{|c|} 
\hline a_1\\ \hline \vdots \\ \hline a_{r}\\ \hline \end{array}\in B^{r,1}$
and $c_2=\begin{array}{|c|} \hline b_1\\ \hline \vdots \\ \hline b_{r}\\ \hline 
\end{array}$, so that we have $a_i\le b_i$ 
for $1\le i \le r$.  Let $(\nu,J)=\Phi(b)$.
We want to show that $(\nu,J) \in \rcls(\RC(L,\la))$. To do that by definition 
of $\rcls$ it is enough to show that $(\nu,J)^{(r)}$  has no singular string of 
length less than $s$. 

Let us introduce some further notation.  Let $\overline{b}=c
\otimes b'$ and $(\nub_0,\Jb_0)=\Phi(c_3\cdots c_{s}\otimes b')$. Define  
$(\nub_i,\Jb_i)=(\rclb^{-1}\circ \delta^{-1})^{i-1}\circ \delta^{-1} (\nub_0,\Jb_0)$ 
for $1\le i \le r$ and let 
$\s_i^{(a)}$ be the length of the singular strings associated to $b_i$. Similarly 
define $(\nu_i,J_i)= (\rclb^{-1}\circ \delta^{-1})^{i-1}\circ \delta^{-1}(\nu_0,J_0)$ 
for $1\le i \le r$ and let  $s_i^{(a)}$ be the length of the singular strings 
associated to $a_i$ where $(\nu_0,J_0)=\Phi(\overline{b})$. The change of 
vacancy number from $(\nub_0,\Jb_0)$ to $(\nub_i,\Jb_i)$ is given by
\begin{equation}\label{eq:1 change in vac}
p_k^{(a)}(\nub_i)= p_k^{(a)}(\nub_0)+\sum_{m=1}^{i}\chi(\s_{m}^{(a-1)}< k \le 
\s_{m}^{(a)})-\sum_{m=1}^{i}\chi( \s_{m}^{(a)}<k \le \s_{m}^{(a+1)}),
\end{equation}
and the change of vacancy number from $(\nub_0,\Jb_0)$ to $(\nu_i,J_i)$ is given by
\begin{equation}\label{eq:2 change in vac}
\begin{split}
p_k^{(a)}(\nu_i)&= p_k^{(a)}(\nub_0)+ \sum_{m=1}^{r}\chi(
 \s_{m}^{(a-1)}< k \le \s_{m}^{(a)})-\sum_{m=1}^{r}\chi( \s_{m}^{(a)}<k \le 
 \s_{m}^{(a+1)})\\ 
 &-\delta_{a,r}\chi (k<s-1)+\sum_{m=1}^{i}\chi(s_{m}^{(a-1)}< k \le 
s_{m}^{(a)})-\sum_{m=1}^{i}\chi( s_{m}^{(a)}<k \le s_{m}^{(a+1)}).
\end{split}
\end{equation}
 
Using this we will show that $s_i^{(a)}>\s_i^{(a)}$ for all $a\ge i$ and $1\le i\le r$ 
by induction on $i$.
Note that by \eqref{eq:1 change in vac}  in $(\nu_0,J_0)^{(a)}$ the strings of
length $\s_i^{(a)}+1$ remain singular for all $i,a$. Since $a_1\le b_1$ we 
have $s_1^{(a)}>\s_1^{(a)}$ for all $a$, this starts the induction. Let  $s_i^{(a)}
>\s_i^{(a)}$ for all $a$ and for $1\le i \le k$. Then by induction hypothesis and 
\eqref{eq:2 change in vac} in $(\nu_k,J_k)^{(a)}$ the strings of length 
$\s_i^{(a)}+1$ remain singular for all $a$ and $k+1\le i \le r$,
which implies that $s_{k+1}^{(a)}\ge \s_{k+1}^{(a)}+1$. Hence $s_{k+1}^{(a)}
> \s_{k+1}^{(a)}$  which proves our claim by induction.
In particular $s_r^{(r)}>\s_r^{(r)}$. By induction $(\nub_r,\Jb_r)^{(r)}$ has no 
singular string of length strictly less than $s-1$, so $\s_r^{(r)}\ge s-1$ which
implies $s_r^{(r)}\ge s$. But note that by construction of the algorithm 
$s_r^{(a)}=0$ for $1\le a\le r-1$ and the change of vacancy numbers
from $(\nu_{r-1},J_{r-1})$ to $(\nu_{r},J_{r})=(\nu,J)$ is given by,
\begin{equation*}
p_k^{(a)}(\nu)=p_k^{(a)}(\nu_{r-1})+\chi(s_{r}^{(a-1)}< k \le s_{r}^{(a)})
                              -\chi( s_{r}^{(a)}<k \le s_{r}^{(a+1)}).
\end{equation*}                              
This implies that $(\nu,J)^{(r)}$ has no singular string less than $s_r^{(r)}$ which
means $(\nu,J)^{(r)}$ has no singular string less than $s$ and we are done.

Conversely let $(\nu,J)\in \rcls(\RC(L,\la))$ and $b=\Phi^{-1}(\nu,J)=c_1\otimes 
c \otimes b'$, same notation as before. We will show that $a_i\le b_i$ for 
$1\le i \le r$. Set $(\nu_i,J_i)=(\delta \circ \lb)^{r-i}(\nu,J)$ for $1\le i\le r$ and 
set  $(\nu_0,J_0)=\delta \circ (\delta \circ \lb)^{r-1}(\nu,J)$.  Let us denote 
the length of the string selected by $\delta$ in $(\nu_i,J_i)^{(a)}$ by $\ell_i^{(a)}$. 
Similarly set $(\nub,\Jb)=\rcls (\nu_0,J_0)$ and  $(\nub_i,\Jb_i)=
(\delta \circ \lb)^{r-i}(\nub,\Jb)$ for $1\le i\le r$ and  $(\nub_0,\Jb_0)=
\delta \circ (\delta \circ \lb)^{r-1}(\nub,\Jb)$. Denote the length of the 
string selected by $\delta$ in $(\nub_i,\Jb_i)^{(a)}$ by $\ellb_i^{(a)}$. 
We claim that $\ell_i^{(a)}> \ellb_i^{(a)}$ for all $1\le i \le r$ and all 
$i\le a \le n$. We will show this by reverse induction on $i$.

First note that the change in vacancy number from $(\nu,J)$ to $(\nu_i,J_i)$ is 
given by
\begin{equation}\label{eq:lrs1 change in vac}
p_k^{(a)}(\nu_i)= p_k^{(a)}(\nu)-\sum_{m=i+1}^{r}\chi(\ell_{m}^{(a-1)}\le k < 
\ell_{m}^{(a)})+\sum_{m=i+1}^{r}\chi( \ell_{m}^{(a)}\le k < \ell_{m}^{(a+1)}).
\end{equation}
The change in vacancy number from $(\nu,J)$ to $(\nub_i,\Jb_i)$ is given by
\begin{equation}\label{eq:lrs2 change in vac}
\begin{split}
p_k^{(a)}&(\nub_i)= p_k^{(a)}(\nu)-\sum_{m=1}^{r}\chi(\ell_{m}^{(a-1)}\le k < 
\ell_{m}^{(a)})+\sum_{m=1}^{r}\chi( \ell_{m}^{(a)}\le k < \ell_{m}^{(a+1)})\\
&+\delta_{a,r}\chi (k<s-1)-\sum_{m=i+1}^{r}\chi(\ellb_{m}^{(a-1)}\le k < 
\ellb_{m}^{(a)})+\sum_{m=i+1}^{r}\chi( \ellb_{m}^{(a)}\le k < \ellb_{m}^{(a+1)}).
\end{split}
\end{equation}
\eqref{eq:lrs1 change in vac} implies that $\ell_i^{(a)}<\ell_{i-1}^{(a)}$ and   
the string of length $\ell_j^{(a)}-1$ remains singular in $(\nu_i,J_i)^{(a)}$ 
for $i+1\le j\le r$. 
Recall that $(\nu,J)^{(r)}$ has no singular string of length less than $s$. So,
$\ell_r^{(r)}\ge s$.  By construction of the algorithm $\ellb_r^{(a)}=1$ for 
$1\le a\le r-1$. By induction $(\nub,\Jb)^{(r)}$ has no singular string of 
length less than $s-1$ and hence by \eqref{eq:lrs2 change in vac} $s-1\le 
\ellb_r^{(r)}< \ell_r^{(r)}$ since the string of length $\ell_r^{(r)}-1\ge s-1$ is 
singular. Now by using \eqref{eq:lrs1 change in vac} the 
algorithm of $\delta$ 
acting on $(\nub,\Jb)$ gives that  $\ellb_r^{(a)}< \ell_r^{(a)}$ for $a\ge r$. This
starts the induction. Suppose  $\ell_i^{(a)}> \ellb_i^{(a)}$ for all $k\le i \le r$ 
and all $i< a \le n$. 
Induction hypothesis along with \eqref{eq:lrs2 change in vac} implies that 
in $(\nub_{k-1},\Jb_{k-1})^{(a)}$ we have
$\ellb_i^{(a)}<\ellb_{i-1}^{(a)}$ for $i\ge k+1$ and   the string of length 
$\ell_j^{(a)}-1$ remains singular  for $1\le j\le k-1$. Therefore 
$\ellb_{k-1}^{(a)}=1$ for $1\le a \le k-2$ and in $(\nub_{k-1},\Jb_{k-1})^{(k-1)}$,
the smallest singular string  we know is of length $\ell_{k-1}^{(k-1)}-1$. Hence 
$\ellb_{k-1}^{(k-1)}\le \ell_{k-1}^{(k-1)}-1< \ell_{k-1}^{(k-1)}$. Then by using 
\eqref{eq:lrs2 change in vac} the algorithm of $\delta$ acting on 
$(\nub_k,\Jb_k)$ gives that  $\ellb_{k-1}^{(a)}< \ell_{k-1}^{(a)}$ for $a>k-1$. 
This proves our claim.      

But  $\ell_i^{(a)}> \ellb_i^{(a)}$ for all 
$1\le i \le r$ and all $i\le a \le n$ implies $a_i\le b_i$. So we are done.

\section{Proof of Theorem~\ref{thm:commute}}\label{appn:crystal}
 
In this section we prove that the crystal operators on paths and rigged
configurations commute with the bijection $\Phi$. A detailed verification of this proof
and its extension to type $D$ is given in~\cite{Sa:2013}.

The following Lemma is a result of~\cite[Lemma 3.11]{KSS:2002} about the convexity of
the vacancy numbers.
\begin{lemma}\label{lem:convexity}({\bf Convexity})
Let $(\nu,J)\in \RC(L)$.
\begin{enumerate}
\item  For all $i,k \ge 1$ we have
 $-p_{k-1}^{(i)}(\nu)+2p_{k}^{(i)}(\nu)-p_{k+1}^{(i)}(\nu)\ge m_{k}^{(i-1)}(\nu)-
 2m_k^{(i)}(\nu)+m_{k}^{(i+1)}(\nu)$.
\item Let $m_k^{(i)}(\nu)=0$ for $a<k<b$. Then
 $p_{k}^{(i)}(\nu)\ge \min(p_{a}^{(i)}(\nu), p_{b}^{(i)}(\nu))$.
\item Let $m_k^{(i)}(\nu)=0$ for $a<k<b$. If $p_a^{(i)}(\nu)=p_{a+1}^{(i)}(\nu)$ and 
$p_{a+1}^{(i)}(\nu)\le p_{b}^{(i)}(\nu)$
then $p_{a+1}^{(i)}(\nu)=p_{k}^{(i)}(\nu)$ for all $a\le k\le b$.
\item Let $m_k^{(i)}(\nu)=0$ for $a<k<b$. If $p_b^{(i)}(\nu)=p_{b-1}^{(i)}(\nu)$ and 
$p_{b-1}^{(i)}(\nu)\le p_{a}^{(i)}(\nu)$
then $p_{b-1}^{(i)}(\nu)=p_{k}^{(i)}(\nu)$ for all $a\le k\le b$.
\end{enumerate}
\end{lemma}
\begin{proof} The proof of (1) is given in \cite[Appendix]{KS:1998} 
(see also~\eqref{eq:p ineq}), (2) follows from repeated use of (1), and the proof of (3)
and (4) follow from (1) and (2).
\end{proof}  

\begin{lemma}\label{lem:delta commute}
Let $B=B^{1,1}\otimes B'$ and let $L$ and $L'$ be the multiplicity 
arrays of $B$ and $B'$. For $1\le i<n$ the following diagrams 
commute if $\ft_i$ is always defined:
\begin{equation} \label{eq:delta commute}
\begin{CD}
\RC(L) @>{\delta}>> \RC(L') \\
@V{\ft_i}VV @VV{\ft_i}V \\
\RC(L) @>>{\delta}> \RC(L')
\end{CD}
\qquad \quad
\begin{CD}
\RC(L) @>{\delta}>> \RC(L') \\
@V{\et_i}VV @VV{\et_i}V \\
\RC(L) @>>{\delta}> \RC(L')
\end{CD}
\end{equation} 
\end{lemma}   

\begin{proof} We prove~\eqref{eq:delta commute} for $\ft_i$ here; the proof for $\et_i$ is 
similar. Let us introduce some notation. Let $(\nu,J)\in \RC(L)$ and let $\ell^{(a)}$ be 
the length of the singular string selected by $\delta$ in $(\nu,J)^{(a)}$ for $1\le a<n$. 
Let $(\nub,\Jb)=\delta(\nu,J)$ and $(\nut,\Jt)=\ft_i(\nu,J)$.  Let $\ellt^{(a)}$ be 
the length of the singular string selected by $\delta$ in $(\nut,\Jt)^{(a)}$ for 
$1\le a<n$ and $\ell$ (respectively $\ellb$) be the length of the string selected by 
$\ft_i$ in $(\nu,J)^{(i)}$ (respectively in $(\nub,\Jb)^{(i)}$). A string of length $k$ 
and label $x_k$ in $(\nu,J)^{(a)}$ is denoted by $(k,x_k)$. 

Using the definition of $\ft_i$ it is easy to see that the diagram \eqref{eq:delta commute}
commutes except when $\ell^{(i-1)}-1\le \ell\le \ell^{(i)}$. We list the 
nontrivial cases as follows:
\begin{enumerate} 
\item[(a)] $\ell^{(i-1)}<\infty, \ell^{(i)}=\infty$, $\ell+1\ge \ell^{(i-1)}$. 
\item[(b)] $\ell^{(i)}<\infty, \ell^{(i-1)}\le \ell+1\le \ell^{(i)}$. 
\item[(c)] $\ell^{(i)}<\infty$ and $\ell^{(i)}=\ell$.
\end{enumerate}
Note that since $\ft_i$ fixes all the colabels, the singular strings (except the new
string of length $\ell+1$) remain singular under the action of $\ft_i$.  Let
$(\ell,x_{\ell})$ be the string selected by $\ft_i$ in $(\nu,J)^{(i)}$.
The new string of length $\ell+1$ can be singular in $(\nut,\Jt)^{(i)}$ only if 
$p_{\ell+1}^{(i)}(\nu)=x_{\ell}+1$. Also note that by the definition of $\ft_i$ if 
$m_k^{(i)}(\nu)>0$ and $(k,x_k)$ is a string in $(\nu,J)^{(i)}$
then
\begin{equation}\label{eq:ft condition}
\begin{split}
x_{\ell} &<x_k\le p_{k}^{(i)}(\nu), \quad \text{if $k>\ell$},\\
x_{\ell} &\le x_k\le p_{k}^{(i)}(\nu),  \quad \text{if $k<\ell$}.
\end{split}
\end{equation} 
Let us now consider the above cases.

\noindent
\textbf{Case (a):} If the new string of length $\ell+1$ in $(\nut,\Jt)^{(i)}$ is 
nonsingular, then \eqref{eq:delta commute}  commutes trivially. Let us 
consider the case when the new string of length $\ell+1$ in 
$(\nut,\Jt)^{(i)}$ is singular. We have  $p_{\ell+1}^{(i)}(\nu)=x_{\ell}+1$ 
and since $\ell^{(i-1)}<\infty, \ell^{(i)}=\infty$ we have $ p_{j}^{(i)}(\nub)
=p_{j}^{(i)}(\nu)-1$  for $j\ge \ell^{(i-1)}$.  
In particular $p_{\ell+1}^{(i)}(\nub)=p_{\ell+1}^{(i)}(\nu)-1=x_{\ell}$.
The labels in $(\nub,\Jb)^{(i)}$ are the same as in $(\nu,J)^{(i)}$. Hence $\ellb=\ell$,  
but the result is not a valid rigged configuration since  
$p_{\ell+1}^{(i)}(\nub)-2<x_{\ell}-1$.
So, $\ft_i(\nub,\Jb)$ is undefined, which contradicts the assumptions of 
Lemma~\ref{lem:delta commute}.

\noindent 
\textbf{Case (b):} If the new string of length $\ell+1$ in 
$(\nut,\Jt)^{(i)}$ is singular, we show that the following conditions hold:
\begin{enumerate}
 \item[(i)] $p_{\ell^{(i)}-1}^{(i)}(\nub)\le x_{\ell}$;
 \item[(ii)] $m_j^{(i+1)}(\nu)=0$ for $\ell<j<\ell^{(i)}$.     
\end{enumerate} 
The above conditions imply that diagram \eqref{eq:delta commute} with $\ft_i$ 
commutes for the following reason. Condition (i) implies that $\ft_i$ acts on the new 
string of length $\ell^{(i)}-1$ in $(\nub,\Jb)^{(i)}$. Condition (ii) implies that  
if $\ell^{(i+1)}<\infty$ then $\ellt^{(i+1)}=\ell^{(i+1)}$. Hence $\ellt^{(a)}=\ell^{(a)}$ 
for $ a \ne i$ and $\ellt^{(i)}=\ell+1$. This gives $\ft_i\circ \delta(\nu,J)=
\delta \circ \ft_i(\nu,J)$. 
 
If  the new string of length $\ell+1$ in $(\nut,\Jt)^{(i)}$ is nonsingular 
then the diagram \eqref{eq:delta commute} with $\ft_i$ commutes if $\ft_i$ acts on 
the same string of length $\ell$ in $(\nub,\Jb)^{(i)}$ as it did on $(\nu,J)^{(i)}$. 
In this case if $(\ell^{(i-1)}-1, p_{\ell^{(i)}-1}^{(i)}(\nub))$ is the new string created 
by $\delta$ we need to show that $x_{\ell}<p_{ \ell^{(i)}-1}^{(i)}(\nub)$.
 
Let us now consider the proof of conditions (i) and (ii) in the case when the new string 
of length $\ell+1$ in $(\nut,\Jt)^{(i)}$ is singular. Note that $p_{\ell+1}^{(i)}(\nu)=
x_{\ell}+1\le x_j$ for $j>\ell$ and $m_j^ {(i)}(\nu)>0$ by~\eqref{eq:ft condition}.
In particular if $m_{\ell+1}^ {(i)}(\nu)>0$ and $(\ell+1,x_{\ell+1})$ is 
a string in $(\nu,J)^{(i)}$  then $p_{\ell+1}^{(i)}(\nu)\le x_{\ell+1}\le 
p_{\ell+1}^{(i)}(\nu)$. This implies $p_{\ell+1}^{(i)}(\nu)=x_{\ell+1}$, 
hence $(\ell+1, x_{\ell+1})$ is a singular string which is a contradiction 
if $\ell^{(i-1)}\le \ell+1< \ell^{(i)}$. If $\ell+1=\ell^{(i)}$, it is easy to
see that~\eqref{eq:delta commute} commutes. Hence we may assume that $\ell+1< \ell^{(i)}$,
so that $m_{\ell+1}^{(i)}(\nu)=0$. 
Let $k>\ell$ be smallest so that $m_{k}^ {(i)}(\nu)>0$. 
Then by Lemma ~\ref{lem:convexity}  (2) we have
\begin{equation} \label{eq:convex1}
p_{\ell+1}^{(i)}(\nu)\ge \min(p_{\ell}^{(i)}(\nu), p_{k}^{(i)}(\nu)).
\end{equation}
If $p_{\ell}^{(i)}(\nu)>p_{k}^{(i)}(\nu)$ then by \eqref{eq:convex1} we 
get $p_{\ell+1}^{(i)}(\nu)\ge p_{k}^{(i)}(\nu)$. But 
\begin{equation}\label{eq:convex2}
\begin{split}
p_{\ell+1}^{(i)}(\nu) &\le x_k< p_{k}^{(i)}(\nu) \quad \text{if $\ell<k<\ell^{(i)}$},\\
p_{\ell+1}^{(i)}(\nu) &\le x_k= p_{k}^{(i)}(\nu) \quad \text{if $k=\ell^{(i)}$}.
\end{split}
\end{equation}
Hence $k=\ell^{(i)}$ which implies $p_{\ell+1}^{(i)}(\nu)=p_{\ell^{(i)}}^{(i)}(\nu)
<p_{\ell}^{(i)}(\nu)$ and $m_{j}^ {(i)}(\nu)=0$ for $\ell<j<\ell^{(i)}$. But now 
using Lemma ~\ref{lem:convexity} (1) we get
the following contradiction:
\begin{equation*}
0>-p_{\ell}^{(i)}(\nu)+2p_{\ell+1}^{(i)}(\nu)-p_{\ell+2}^{(i)}(\nu)\ge 
m_{\ell+1}^{(i-1)}(\nu)+ m_{\ell+1}^{(i+1)}(\nu)\ge 0.
\end{equation*}
Hence $p_{\ell}^{(i)}(\nu)\le p_{k}^{(i)}(\nu)$ and by \eqref{eq:convex1} 
we $p_{\ell+1}^{(i)}(\nu)\ge p_{\ell}^{(i)}(\nu)$. Recall that we have 
$m_{\ell+1}^{(i)}(\nu)=0$ and 
\begin{align*}
p_{\ell+1}^{(i)}(\nu)&=x_{\ell}+1\le p_{\ell}^{(i)}(\nu),
 &&\text{if $\ell^{(i-1)}<\ell<\ell^{(i)}$ or 
         $(\ell, x_{\ell})$ is nonsingular},\\
p_{\ell+1}^{(i)}(\nu)&=x_{\ell}+1=p_{\ell}^{(i)}(\nu)+1,
 &&\text{if $\ell=\ell^{(i-1)}-1$ and $(\ell, x_{\ell})$ is singular}.
\end{align*}
This gives us two possible situations:
\begin{enumerate}
\item $p_{\ell+1}^{(i)}(\nu)=p_{\ell}^{(i)}(\nu)$ if $\ell^{(i-1)}<\ell<\ell^{(i)}$ or 
                            $(\ell, x_{\ell})$ is nonsingular,
\item  $p_{\ell+1}^{(i)}(\nu)= p_{\ell}^{(i)}(\nu)+1$  if $\ell=\ell^{(i-1)}-1$ and 
                          $(\ell, x_{\ell})$ is singular.
\end{enumerate}

In situation (1) using Lemma~\ref{lem:convexity} (3) we get $p_{\ell+1}^{(i)}(\nu)=
p_{j}^{(i)}(\nu)$ for $\ell+1< j \le k$. Using \eqref{eq:convex2} this implies 
$k=\ell^{(i)}$ and by convexity we get condition (ii). Also this
gives $p_{\ell^{(i)}-1}^{(i)}(\nu)=p_{\ell+1}^{(i)}(\nu)=x_{\ell}+1$ and hence 
$p_{\ell^{(i)}-1}^{(i)}(\nub)=x_{\ell}$, which proves condition (i).                        

In situation (2)\footnote{We thank Reiho Sakamoto for pointing out a typo in a previous version of this paragraph.
See also~\cite[Section 5.6.3]{Sa:2013}.},
since $m_{\ell+1}^{(i-1)}(\nu)>0$ and $m_{\ell+1}^{(i)}(\nu)=0$, by the convexity
of the vacancy numbers and~\eqref{eq:convex2} we find that 
$p_{\ell+1}^{(i)}(\nu)=p_{\ell+2}^{(i)}(\nu)=\cdots=p_k^{(i)}(\nu)$ with $k=\ell^{(i)}$.
This implies in particular that $p_{\ell^{(i)}-1}^{(i)}(\nu)=p_{\ell+1}^{(i)}(\nu)
=p_\ell^{(i)}(\nu)+1=x_\ell+1$ since $(\ell,x_\ell)$ is singular.
By assumption $\ell^{(i-1)}=\ell+1<\ell^{(i)}$, so that $p_{\ell^{(i)}-1}^{(i)}(\nub)
=p_{\ell^{(i)}-1}(\nu)-1=x_\ell$, which proves condition (i). By the same arguments as in
situation (1), condition (ii) follows.

Now let us consider the case when  the new string of length $\ell+1$ in $(\nut,\Jt)^{(i)}$ is 
nonsingular. If $\ell+1=\ell^{(i)}$ the commutation of~\eqref{eq:delta commute} is again
fairly easy to see. Hence assume that $\ell+1<\ell^{(i)}$. Then we have 
$p_{\ell^{(i)}-1}^{(i)}(\nub)=p_{\ell^{(i)}-1}^{(i)}(\nu)-1$. If 
$m_{\ell^{(i)}-1}^{(i)}(\nu)>0$  and $(\ell^{(i)}-1,x_{\ell^{(i)}-1})$
is a string in $(\nu,J)^{(i)}$  then $x_{\ell^{(i)}-1}<p_{\ell^{(i)}-1}^{(i)}(\nu)$ since 
$\ell^{(i-1)}\le \ell+1\le \ell^{(i)}-1<\ell^{(i)}$. Hence by \eqref{eq:ft condition} we 
have $x_{\ell}<x_{\ell^{(i)}-1}<p_{\ell^{(i)}-1}^{(i)}(\nu)$ which implies $x_{\ell}<
p_{\ell^{(i)}-1}^{(i)}(\nub)$ and we are done.

If $m_{\ell^{(i)}-1}^{(i)}(\nu)=0$ let $\ell \le j<\ell^{(i)}-1$ be smallest such that 
$m_j^{(i)}(\nu)>0$.
By Lemma ~\ref{lem:convexity} (2) we get 
\begin{equation}\label{eq:min}
p_{\ell^{(i)}-1}^{(i)}(\nu)\ge \min( p_{j}^{(i)}(\nu), p_{\ell^{(i)}}^{(i)}(\nu)).
\end{equation}
Note that if $\ell<j<\ell^{(i)}$ then the string $(j,x_j)$ in $(\nu,J)^{(i)}$ is 
nonsingular and therefore  $p_{j}^{(i)}(\nu)>x_j>x_{\ell}$ by \eqref{eq:ft condition}.
Also if $(\ell^{(i)},x_{\ell^{(i)}})$ is the singular string $p_{\ell^{(i)}}^{(i)}(\nu)=
x_{\ell^{(i)}}>x_{\ell}$ by  \eqref{eq:ft condition}.
So $\min( p_{j}^{(i)}(\nu), p_{\ell^{(i)}}^{(i)}(\nu))\ge x_{\ell}+1$. 
Hence by \eqref{eq:min} $p_{\ell^{(i)}-1}^{(i)}(\nu)\ge x_{\ell}+1$.
Suppose $p_{\ell^{(i)}-1}^{(i)}(\nu)= x_{\ell}+1$. Since $p_{j}^{(i)}(\nu)>x_{\ell}+1$
we get by \eqref{eq:min}  $x_{\ell}+1=p_{\ell^{(i)}-1}^{(i)}(\nu)\ge 
p_{\ell^{(i)}}^{(i)}(\nu) \ge x_{\ell}+1$
which implies $p_{\ell^{(i)}-1}^{(i)}(\nu)= p_{\ell^{(i)}}^{(i)}(\nu)$. Since 
$p_{\ell^{(i)}-1}^{(i)}(\nu)=x_{\ell}+1\le  p_{a}^{(i)}(\nu)$ for all $j<a<\ell^{(i)}$ by
Lemma ~\ref{lem:convexity} (4) we get $p_{j}^{(i)}(\nu)=x_{\ell}+1$ which is a
contradiction. Hence $p_{\ell^{(i)}-1}^{(i)}(\nu)> x_{\ell}+1$ and we get $x_{\ell}<
p_{\ell^{(i)}-1}^{(i)}(\nub)$ as desired.

Let us consider the case $j=\ell$. If the string $(\ell, x_{\ell})$ is nonsingular
by similar argument as in the previous case we have that $p_{\ell^{(i)}-1}^{(i)}(\nu)\ge
x_{\ell}+1$. Suppose   $p_{\ell^{(i)}-1}^{(i)}(\nu)=
x_{\ell}+1$. By \eqref{eq:min} if $p_{\ell^{(i)}-1}^{(i)}(\nu)\ge p_{\ell^{(i)}}^{(i)}(\nu)\ge 
x_{\ell}+1$ we get as before that $p_{\ell^{(i)}-1}^{(i)}(\nu)= p_{\ell^{(i)}}^{(i)}(\nu)$. 
Using Lemma~\ref{lem:convexity} (4)  we can show as before that  $p_{\ell+1}^{(i)}(\nu)
=x_{\ell}+1$ which is a contradiction since the string of length $\ell+1$ is not singular
in $(\nut,\Jt)^{(i)}$. By \eqref{eq:min} if $p_{\ell^{(i)}-1}^{(i)}(\nu)\ge 
p_{\ell}^{(i)}(\nu)\ge x_{\ell}+1$ we get  $p_{\ell^{(i)}-1}^{(i)}(\nu)= p_{\ell}^{(i)}(\nu)
=x_{\ell}+1$. This implies that $p_{\ell^{(i)}-1}^{(i)}(\nu)\le p_{a}^{(i)}(\nu)$ for all 
$a>\ell$. If we use this in Lemma~\ref{lem:convexity} (1) for $k=\ell^{(i)}-1$ we get
$p_{\ell^{(i)}-1}^{(i)}(\nu)=p_{\ell^{(i)}}^{(i)}(\nu)$ and then using 
Lemma~\ref{lem:convexity} (4) we get $p_{\ell+1}^{(i)}(\nu)=x_{\ell}+1$ which is a 
contradiction as before.

Hence the only case left to be considered is when $j=\ell=\ell^{(i-1)}-1$ and the string 
$(\ell, x_{\ell})$ is singular in $(\nu,J)^{(i)}$. Here $ \min( p_{\ell}^{(i)}(\nu), 
p_{\ell^{(i)}}^{(i)}(\nu))= p_{\ell}^{(i)}(\nu)$ and therefore by \eqref{eq:min}
$p_{\ell^{(i)}-1}^{(i)}(\nu)\ge x_{\ell}$. Suppose $p_{\ell^{(i)}-1}^{(i)}(\nu)=x_{\ell}$.
Since $p_{\ell^{(i)}}^{(i)}(\nu)\ge x_{\ell}+1$ we have $p_{\ell^{(i)}-1}^{(i)}(\nu)<
p_{\ell^{(i)}}^{(i)}(\nu)$. Also, $p_{\ell^{(i)}-1}^{(i)}(\nu)\ge \min(p_{\ell}^{(i)}(\nu),
p_{\ell^{(i)}}^{(i)}(\nu))=p_{\ell}^{(i)}(\nu)=x_\ell=p_{\ell^{(i)}-1}^{(i)}(\nu)$. Using 
this in Lemma~\ref{lem:convexity} (1) for $k=\ell^{(i)}-1$ we get the following contradiction:
\begin{equation}\label{eq:contradiction1} 
 0>-p_{\ell^{(i)}-2}^{(i)}(\nu)+2p_{\ell^{(i)}-1}^{(i)}(\nu)-p_{\ell^{(i)}}^{(i)}(\nu)\ge 
m_{\ell^{(i)}-1}^{(i-1)}(\nu)+ m_{\ell^{(i)}-1}^{(i+1)}(\nu)\ge 0.
\end{equation}
Hence $p_{\ell^{(i)}-1}^{(i)}(\nu)> x_{\ell}$. Suppose $p_{\ell^{(i)}-1}^{(i)}(\nu)=
 x_{\ell}+1$. Here $p_{\ell^{(i)}-1}^{(i)}(\nu)\le p_{\ell^{(i)}}^{(i)}(\nu)$. If 
 $p_{\ell^{(i)}-1}^{(i)}(\nu)= p_{\ell^{(i)}}^{(i)}(\nu)$ as before we can show that
$p_{\ell+1}^{(i)}(\nu)=x_{\ell}+1$, which is a contradiction. Suppose
$p_{\ell^{(i)}-1}^{(i)}(\nu)< p_{\ell^{(i)}}^{(i)}(\nu)$ then $p_{\ell^{(i)}-2}^{(i)}(\nu)
\ge \min( p_{\ell}^{(i)}(\nu),p_{\ell^{(i)}}^{(i)}(\nu))=p_{\ell}^{(i)}(\nu)=x_{\ell}=
p_{\ell^{(i)}-1}^{(i)}(\nu)-1$. If $p_{\ell^{(i)}-2}^{(i)}(\nu)>p_{\ell^{(i)}-1}^{(i)}(\nu)$
we again get  the contradiction  \eqref{eq:contradiction1}. If 
$p_{\ell^{(i)}-2}^{(i)}(\nu)=p_{\ell^{(i)}-1}^{(i)}(\nu)$ using 
Lemma~\ref{lem:convexity} (1) for $k=\ell^{(i)}-1$ we get  $p_{\ell^{(i)}}^{(i)}(\nu)
=p_{\ell^{(i)}-1}^{(i)}(\nu)$ which is a contradiction to our assumption. Hence 
$p_{\ell^{(i)}-1}^{(i)}(\nu)> x_{\ell}+1$ giving $x_{\ell}<p_{\ell^{(i)}-1}^{(i)}(\nub)$.

\noindent
\textbf{Case (c):}  Note that since $\ft_i$ acts on the string   $(\ell, x_{\ell})$ in 
$(\nu,J)^{(i)}$ we have  
\begin{equation}\label{eq:convex3}
p_{\ell+1}^{(i)}(\nu)\ge x_{\ell}+1=p_{\ell}^{(i)}(\nu)+1.
\end{equation} 
If $\ft_i$ and $\delta$ select the same string of length $\ell$ in $(\nu,J)^{(i)}$ 
then $m_{\ell}^{(i)}(\nu)= 1$.
But if $\ft_i$ and $\delta$ select different strings of length $\ell$ in $(\nu,J)^{(i)}$ 
then  $m_{\ell}^{(i)}(\nu)>1$. We will consider each of these two cases separately.

If  $m_{\ell}^{(i)}(\nu)>1$  let $(\ell, x_{\ell})$ be the string selected by $\ft_i$ 
and $(\ell, p_{\ell}^{(i)}(\nu))$ be the string selected by $\delta$ in $(\nu,J)^{(i)}$. 
Note that $x_{\ell}\le p_{\ell}^{(i)}(\nu)$. To prove that the 
diagram \eqref{eq:delta commute} with $\ft_i$ commutes it is enough to show that 
$\ft_i$ acts on the same string $(\ell, x_{\ell})$ in $(\nub,\Jb)^{(i)}$ as it did 
in $(\nu,J)^{(i)}$. Hence it suffices to show that the new label in $(\nub,\Jb)^{(i)}$ 
satisfies $p_{\ell-1}^{(i)}(\nub)\ge x_{\ell}$. Note that
\begin{align*}
p_{\ell-1}^{(i)}(\nub)&=p_{\ell-1}^{(i)}(\nu)-1 &&\text{if $\ell>\ell^{(i-1)}$,}\\
p_{\ell-1}^{(i)}(\nub)&=p_{\ell-1}^{(i)}(\nu) &&\text{if $\ell=\ell^{(i-1)}$}.
\end{align*}

If $m_{\ell-1}^{(i)}(\nu)>0$  let $(\ell-1,x_{\ell-1})$ be a string in $(\nu,J)^{(i)}$. 
Then
\begin{align*}
x_{\ell}&\le x_{\ell-1}< p_{\ell-1}^{(i)}(\nu) &&\text{if $\ell>\ell^{(i-1)}$,}\\
x_{\ell}&\le x_{\ell-1}\le p_{\ell-1}^{(i)}(\nu) &&\text{if $\ell=\ell^{(i-1)}$,}
\end{align*}
which implies $p_{\ell-1}^{(i)}(\nub)\ge x_{\ell}$. 

If $m_{\ell-1}^{(i)}(\nu)=0$ let $j<\ell-1$ be largest such that $m_j^{(i)}(\nu)>0$ 
and $(j,x_j)$ be a string in $(\nu,J)^{(i)}$. Then by Lemma~\ref{lem:convexity} (2) we 
have $p_{\ell-1}^{(i)}(\nu)\ge \min(p_{j}^{(i)}(\nu),
p_{\ell}^{(i)}(\nu))$. 

If $p_{j}^{(i)}(\nu)\le p_{\ell}^{(i)}(\nu)$ then using 
\eqref{eq:ft condition} we have
\begin{align*}
p_{\ell-1}^{(i)}(\nu)&\ge p_{j}^{(i)}(\nu)>x_j\ge x_{\ell}
 &&\text{if $\ell^{(i-1)}\le j<\ell-1$,}\\
p_{\ell-1}^{(i)}(\nu)&\ge p_{j}^{(i)}(\nu)\ge x_j\ge x_{\ell}
 &&\text{if $j<\ell^{(i-1)}$.}
\end{align*}
Hence $p_{\ell-1}^{(i)}(\nub)\ge x_{\ell}$ unless 
\begin{equation}\label{eq:special case}
p_{\ell-1}^{(i)}(\nu)=p_{j}^{(i)}(\nu)=x_j= x_{\ell} \le p_{\ell}^{(i)}(\nu) \qquad 
\text{with $j<\ell^{(i-1)} \le \ell-1$.}
\end{equation}
But if this happens by Lemma~\ref{lem:convexity} we get $p_{\ell^{(i-1)}}^{(i)}(\nu)=
 p_{\ell^{(i-1)}-1}^{(i)}(\nu)\le p_{\ell^{(i-1)}+1}^{(i)}(\nu)$. Note that here 
 $m_{\ell^{(i-1)}}^{(i)}(\nu)=0$ and $m_{\ell^{(i-1)}}^{(i-1)}(\nu)>0$. Using all these 
we get the following contradiction:
\begin{equation*}
0\ge -p_{\ell^{(i-1)}-1}^{(i)}(\nu)+2p_{\ell^{(i-1)}}^{(i)}(\nu)
 -p_{\ell^{(i-1)}+1}^{(i)}(\nu)
 \ge m_{\ell^{(i-1)}}^{(i-1)}(\nu)+m_{\ell^{(i-1)}}^{(i-1)}(\nu)\ge 1.
\end{equation*}
This shows that \eqref{eq:special case} can not happen.
 
If $p_{j}^{(i)}(\nu)> p_{\ell}^{(i)}(\nu)$ then $p_{\ell-1}^{(i)}(\nu)\ge 
\min(p_{j}^{(i)}(\nu), p_{\ell}^{(i)}(\nu))=p_{\ell}^{(i)}(\nu)\ge x_{\ell}$. 
Again $p_{\ell-1}^{(i)}(\nub)\ge x_{\ell}$ unless 
\begin{equation}\label{eq:special case2}
p_{\ell-1}^{(i)}(\nu)=p_{\ell}^{(i)}(\nu)= x_{\ell} \qquad 
\text{with $\ell^{(i-1)} \le \ell-1$.}
\end{equation}
But this implies by  Lemma~\ref{lem:convexity}  that $p_{j}^{(i)}(\nu)=p_{\ell}^{(i)}(\nu)
= x_{\ell}$ which is a contradiction to our assumption. Hence \eqref{eq:special case2}
does not occur. This completes the proof when $m_{\ell}^{(i)}(\nu)>1$. 
  
If  $m_{\ell}^{(i)}(\nu)= 1$  we claim that 
\begin{enumerate}
\item[(i)]  $p_{\ell+1}^{(i)}(\nu)=x_{\ell}+1=p_{\ell}^{(i)}(\nu)+1$,
\item[(ii)]  $p_{\ell-1}^{(i)}(\nub)=x_{\ell}$,
\item[(iii)] If $\ell^{(i+1)}<\infty$ then  $\ell+1\le \ell^{(i+1)}$.  
\end{enumerate}
It is easy to see that diagram \eqref{eq:delta commute} with $\ft_i$ commutes if 
our claim is true. Condition (i) implies that the new string 
$(\ell+1,x_{\ell}-1)$ in $(\nut,\Jt)^{(i)}$ is singular and $\ellt^{(i)}=\ell+1$. 
Condition (iii) implies that $\ellt^{(i+1)}=\ell^{(i+1)}$. On the other hand condition
(ii) implies $\ellb=\ell-1$, the new string created by $\delta$ in $(\nub,\Jb)^{(i)}$. 

Let us prove our claims now. Using Lemma~\ref{lem:convexity} (1) we have
\begin{equation*}
(p_{\ell}^{(i)}(\nu)-p_{\ell-1}^{(i)}(\nu))+(p_{\ell}^{(i)}(\nu)-p_{\ell+1}^{(i)}(\nu))\ge 
m_{\ell}^{(i-1)}(\nu)-2+m_{\ell}^{(i+1)}(\nu).
\end{equation*}
which can be rewritten as
\begin{equation}\label{eq:convex4}
(p_{\ell}^{(i)}(\nu)+1-p_{\ell-1}^{(i)}(\nu))+(p_{\ell}^{(i)}(\nu)+1
 -p_{\ell+1}^{(i)}(\nu))\ge m_{\ell}^{(i-1)}(\nu)+m_{\ell}^{(i+1)}(\nu)\ge 0.
\end{equation}
Suppose  $\ell^{(i-1)}<\ell=\ell^{(i)}$. If $m_{\ell-1}^{(i)}(\nu)>0$ then the string 
$(\ell-1,x_{\ell-1})$ is nonsingular and hence by \eqref{eq:ft condition} 
$p_{\ell}^{(i)}(\nu)=x_{\ell}\le x_{\ell-1}< p_{\ell-1}^{(i)}(\nu)$.  If 
$m_{\ell-1}^{(i)}(\nu)=0$  let $j<\ell-1$ be largest such that
$m_j^{(i)}(\nu)>0$. Note that $p_{j}^{(i)}(\nu)\ge x_{\ell}=p_{\ell}^{(i)}(\nu)$, so by 
Lemma~\ref{lem:convexity} (2) we have
$p_{\ell-1}^{(i)}(\nu)\ge \min(p_{j}^{(i)}(\nu),p_{\ell}^{(i)}(\nu))=p_{\ell}^{(i)}(\nu)$.
Hence $p_{\ell-1}^{(i)}(\nu)>p_{\ell}^{(i)}(\nu)$ unless
\begin{equation}\label{eq:special case3}
p_{\ell-1}^{(i)}(\nu)=p_{\ell}^{(i)}(\nu)=p_{j}^{(i)}(\nu)=x_{\ell} 
 \text{with $j<\ell^{(i-1)}<\ell$.}
\end{equation}
But if this happens by Lemma~\ref{lem:convexity} we get $p_{\ell^{(i-1)}}^{(i)}(\nu)
=p_{\ell^{(i-1)}-1}^{(i)}(\nu)=
p_{\ell^{(i-1)}+1}^{(i)}(\nu)$ which gives us the following contradiction since 
$m_{\ell^{(i-1)}}^{(i-1)}(\nu)>0$:
\begin{equation*}
0\ge -p_{\ell^{(i-1)}-1}^{(i)}(\nu)+2p_{\ell^{(i-1)}}^{(i)}(\nu)
-p_{\ell^{(i-1)}+1}^{(i)}(\nu)\ge m_{\ell^{(i-1)}}^{(i-1)}(\nu)
+m_{\ell^{(i-1)}}^{(i-1)}(\nu)\ge 1.
\end{equation*}
Hence  \eqref{eq:special case3} cannot happen and we have $p_{\ell-1}^{(i)}(\nu)>
p_{\ell}^{(i)}(\nu)$. Now using this and \eqref{eq:convex3} in~\eqref{eq:convex4}
we get
\begin{equation*}
0\ge (p_{\ell}^{(i)}(\nu)+1-p_{\ell-1}^{(i)}(\nu))+(p_{\ell}^{(i)}(\nu)+1
 -p_{\ell+1}^{(i)}(\nu))\ge m_{\ell}^{(i-1)}(\nu)+m_{\ell}^{(i+1)}(\nu)\ge 0,
\end{equation*}
which implies $p_{\ell}^{(i)}(\nu)=p_{\ell-1}^{(i)}(\nu)-1$, $p_{\ell+1}^{(i)}(\nu)
=p_{\ell}^{(i)}(\nu)+1$, 
$m_{\ell}^{(i-1)}(\nu)=0$ and $m_{\ell}^{(i+1)}(\nu)=0$. This proves (i) and (iii). Also 
$p_{\ell}^{(i)}(\nu)=p_{\ell-1}^{(i)}(\nu)-1$ implies $p_{\ell-1}^{(i)}(\nub)
=p_{\ell-1}^{(i)}(\nu)-1=p_{\ell}^{(i)}(\nu)=x_{\ell}$.  This proves (ii).

Suppose $\ell^{(i-1)}=\ell=\ell^{(i)}$. This means $m_{\ell}^{(i-1)}(\nu)\ge 1$ and 
as before if $m_{\ell-1}^{(i)}(\nu)>0$ we have $p_{\ell}^{(i)}(\nu)=
x_{\ell}\le x_{\ell-1} \le p_{\ell-1}^{(i)}(\nu)$.  If $m_{\ell-1}^{(i)}(\nu)=0$ 
again as in the previous case we have 
$p_{\ell-1}^{(i)}(\nu)\ge \min(p_{j}^{(i)}(\nu),p_{\ell}^{(i)}(\nu))=p_{\ell}^{(i)}(\nu)$. 
Using this and~\eqref{eq:convex3} in~\eqref{eq:convex4} we get 
$p_{\ell}^{(i)}(\nu)=p_{\ell-1}^{(i)}(\nu)$, $p_{\ell+1}^{(i)}(\nu)=
p_{\ell}^{(i)}(\nu)+1$, $m_{\ell}^{(i-1)}(\nu)=1$ and $m_{\ell}^{(i+1)}(\nu)=0$. 
Note that since $\ell^{(i-1)}=\ell$, $p_{\ell-1}^{(i)}(\nub)=p_{\ell-1}^{(i)}(\nu)
=p_{\ell}^{(i)}(\nu)=x_{\ell}$. So we proved (i), (ii) and (iii).
\end{proof}

\begin{lemma}\label{lem:j commute}
Let $B=B^{r,1}\otimes B'$, $r\ge 2$ and let $L$ be the multiplicity
array of $B$. For $1\le i<n$ the following diagrams commute:
\begin{equation}\label{eq:j commute}
\begin{CD}
\RC(L) @>{\rclb}>> \RC(\lb(L)) \\
@V{\ft_i}VV @VV{\ft_i}V \\
\RC(L) @>>{\rclb}> \RC(\lb(L))
\end{CD}
\qquad \quad
\begin{CD}
\RC(L) @>{\rclb}>> \RC(\lb(L)) \\
@V{\et_i}VV @VV{\et_i}V \\
\RC(L) @>>{\rclb}> \RC(\lb(L))
\end{CD}
\end{equation} 
\end{lemma}

\begin{proof}
Note that if $i>r-1$ then the proof of  \eqref{eq:j commute} is trivial. Suppose 
$1\le i\le r-1$. The proof for $\et_i$ is very similar to the proof for $\ft_i$, 
so here we only prove~\eqref{eq:j commute} for $\ft_i$.
Let $(\nu,J)\in \RC(L)$. Let $(\ell, x_{\ell})$ be the string 
selected by $\ft_i$ in $(\nu,J)^{(i)}$. Let $(\nub,\Jb)=\rclb(\nu,J)$. By 
definition of $\rclb$ we get $(\nub,\Jb)^{(k)}$ by adding a singular string 
of length one to $(\nu,J)^{(k)}$ for $1\le k \le r-1$. Hence to show that 
the diagram~\eqref{eq:j commute} commutes it suffices to show that the label for 
the new singular string of length one in $(\nub,\Jb)^{(i)}$ satisfies
$p_1^{(i)}(\nub)\ge x_{\ell}$. 
Note that $p_1^{(i)}(\nub)= p_1^{(i)}(\nu)$ for all $1\le i\le r-1$.
					    		    	 
If $m_1^{(i)}(\nu)>0$ then $ x_1^{(i)}\ge x_{\ell}$ by 
\eqref{eq:ft condition}.  So,  $p_1^{(i)}(\nub)=
p_1^{(i)}(\nu)\ge x_{1}^{(i)}\ge x_{\ell}$.  
If $m_1^{(i)}(\nu)=0$ let $j$ be smallest such that $m_j^{(i)}(\nu)>0$ and
$(j,x_j)$ be a string in $(\nu,J)^{(i)}$. By Lemma~\ref{lem:convexity} (2) we get
$p_{1}^{(i)}(\nu)\ge \min(p_{0}^{(i)}(\nu), p_{j}^{(i)}(\nu))$. Recall that
$p_{0}^{(i)}(\nu)=0$ and $x_{\ell}\le 0$ by the definition 
of $\ft_i$. So, if $p_{j}^{(i)}(\nu)\ge 0$ then
$p_{j}^{(i)}(\nub)=p_{1}^{(i)}(\nu)\ge 0\ge x_{\ell}$.  If $p_{j}^{(i)}(\nu)< 0$ then  
$p_{1}^{(i)}(\nu)\ge p_{j}^{(i)}(\nu)$. But $p_{j}^{(i)}(\nu)\ge x_{j}
\ge x_{\ell}$. Hence $p_1^{(i)}(\nub)= p_1^{(i)}(\nu)\ge p_{j}^{(i)}(\nu) \ge 
x_{\ell}$ and we are done.
\end{proof}   

\begin{lemma}\label{lem:i commute}
Let $B=B^{r,s}\otimes B'$, $r\ge 1, s\ge 2$ and let $L$ be the multiplicity
array of $B$. For $1\le i<n$ the following diagrams commute:
\begin{equation}\label{eq:i commute}
\begin{CD}
\RC(L) @>{\rcls}>> \RC(\ls(L)) \\
@V{\ft_i}VV @VV{\ft_i}V \\
\RC(L) @>>{\rcls}> \RC(\ls(L))
\end{CD}
\qquad \quad 
\begin{CD}
\RC(L) @>{\rcls}>> \RC(\ls(L)) \\
@V{\et_i}VV @VV{\et_i}V \\
\RC(L) @>>{\rcls}> \RC(\ls(L))
\end{CD}
\end{equation} 
\end{lemma}  
\begin{proof}
Let $(\nu,J)\in \RC(L)$. By definition $\rcls$ only changes the vacancy
numbers in $(\nu,J)^{(r)}$. Hence the proof of this lemma is trivial.  
\end{proof} 

Now we will prove Theorem~\ref{thm:commute}.
\begin{proof}[Proof of Theorem~\ref{thm:commute}]
To prove this theorem we will use a diagram of the form
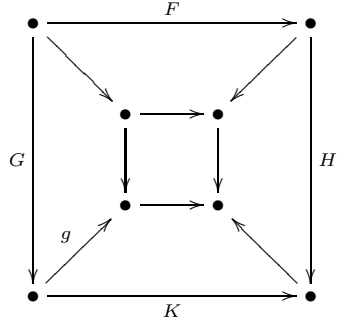
\begin{equation*}
\xymatrix{
 {\bullet} \ar[rrr]^{F} \ar[ddd]_{G} \ar[dr] & & &
        {\bullet} \ar[ddd]^{H} \ar[dl] \\
 & {\bullet} \ar[r] \ar[d] & {\bullet} \ar[d] & \\
 & {\bullet} \ar[r]   & {\bullet}  & \\
 {\bullet} \ar[rrr]_{K} \ar[ur]^{g} & & & {\bullet} \ar[ul]
}
\end{equation*}
We view this diagram as a cube with front face given by the large square. 
By~\cite[Lemma 5.3]{KSS:2002} if the squares given by all the faces of
the cube except the front commute and the map $g$ is injective then the
front face also commutes.

We will prove Theorem~\ref{thm:commute} by using induction on $B$
as we did in the proof of the bijection of Proposition~\ref{prop:bij}. 
First let $B=B^{1,1}\otimes B'$. We prove Theorem~\ref{thm:commute} for $\ft_i$ 
by using Lemma~\ref{lem:delta commute} and the following diagram when 
$f_i$ and $\ft_i$ are defined: 
\begin{equation*}
\xymatrix{
 {\Path(B)} \ar[rrr]^{\Phi} \ar[ddd]_{f_i} \ar[dr] ^{\lh}& & &
        {\RC(L)} \ar[ddd]^{\ft_i} \ar[dl]_{\delta} \\
 & {\Path(B')} \ar[r]^{\Phi} \ar[d]_{f_i} & {\RC(L')} \ar[d] ^{\ft_i}& \\
 & {\Path(B')} \ar[r]^{\Phi}  & {\RC(L')}  & \\
 {\Path(B)} \ar[rrr]_{\Phi} \ar[ur] ^{\lh}& & & {\RC(L)} \ar[ul]_{\delta}
}
\end{equation*}
Note the top and the bottom faces commute by Definition~\ref{def:bij} (1).
The right face commutes by Lemma~\ref{lem:delta commute}. The left face commutes
by definition of $f_i$ on the paths and we know $\lh$ is injective. By induction
hypothesis the back face commutes. Hence the front face must commute.

Let us now prove Theorem~\ref{thm:commute} when not all $f_i$ (resp. $\ft_i$)
in the above diagram are defined. Let $(\nu,J)\in\RC(L)$, $(\nub,\Jb)=\delta(\nu,J)$,
$b=\Phi^{-1}(\nu,J)$ and $b'=\Phi^{-1}(\nub,\Jb)$.
We need to show the following cases:
\begin{enumerate}
\item $f_i(b)$ is defined and $f_i(b')$ is undefined if and only if
      $\ft_i(\nu,J)$ is defined and $\ft_i(\nub,\Jb)$ is undefined.
      In addition $\Phi(f_i(b))=\ft_i(\nu,J)$.
\item $f_i(b)$ is undefined and $f_i(b')$ is defined if and only if
      $\ft_i(\nu,J)$ is undefined and $\ft_i(\nub,\Jb)$ is defined.
\item $f_i(b)$ and $f_i(b')$ are both undefined if and only if
      $\ft_i(\nu,J)$ and $\ft_i(\nub,\Jb)$ are both undefined.
\end{enumerate}

For Case (1) suppose that $\ft_i(\nu,J)=(\nut,\Jt)$ is defined, but $\ft_i(\nub,\Jb)$
is undefined. Then we are in the situation described in Case (a) of 
Lemma~\ref{lem:delta commute}.
That is  $\ell^{(i-1)}<\infty$, $\ell^{(i)}=\infty$, $\ell+1\ge \ell^{(i-1)}$ and 
the new string of length $\ell+1$ is singular in $(\nut,\Jt)^{(i)}$. 
In this situation note that $m_{\ell+1}^{(i)}(\nub)=0$, else $p_{\ell+1}^{(i)}(\nub)\ge x_{\ell+1}>x_{\ell}$ 
by \eqref{eq:ft condition}, which is a contradiction to $p_{\ell+1}^{(i)}(\nub)=x_{\ell}$
as discussed in Case (a) of Lemma~\ref{lem:delta commute}. 
Suppose $j>\ell$ be smallest such that $m_{j}^{(i)}(\nub)>0$. Then
\begin{equation}\label{cond1}
p_j^{(i)}(\nub)\ge x_j>x_{\ell}=p_{\ell+1}^{(i)}(\nub).
\end{equation}
By Lemma~\ref{lem:convexity} (2), $p_{\ell+1}^{(i)}(\nub)\ge \min(p_{\ell}^{(i)}(\nub),
p_j^{(i)}(\nub))$. By \eqref{cond1} this implies $p_{\ell+1}^{(i)}(\nub)\ge 
p_{\ell}^{(i)}(\nub)$.
But $x_{\ell}=p_{\ell+1}^{(i)}(\nub)\ge p_{\ell}^{(i)}(\nub)\ge x_{\ell}$, hence we get
$p_{\ell+1}^{(i)}(\nub)=p_{\ell}^{(i)}(\nub)$. Again by Lemma~\ref{lem:convexity} (3) 
since $m_k^{(i)}(\nub)=0$ for $\ell< k < j$ we get $p_{\ell+1}^{(i)}(\nub)=p_{j}^{(i)}(\nub)$ 
which contradicts \eqref{cond1}.
Hence $m_j^{(i)}(\nub)=0$ for $j>\ell$. Also by Lemma~\ref{lem:convexity} (1) 
$p_{\ell+1}^{(i)}(\nub)=p_{\ell}^{(i)}(\nub)$ with $m_j^{(i)}(\nub)=0$ for $j>\ell$ 
implies that $m_j^{(i+1)}(\nub)=0$ for $j>\ell$. Since  $\nub^{(i+1)}$ and $\nut^{(i+1)}$  
have the same shape we get $m_j^{(i+1)}(\nut)=0$ for $j>\ell$. Hence $\ellt^{(a)}=\ell^{(a)}$ 
for $1\le a \le i-1$, $\ellt^{(i)}=\ell+1$ and $\ellt^{(i+1)}=\infty$.
Therefore we proved that if $\Phi^{-1}(\nub,\Jb)=b' \in B'$ then $\Phi^{-1}(\nu,J)=i 
\otimes b'$ and $\Phi^{-1}(\nut,\Jt)=i+1 \otimes b'$. But $\ft_i(\nub,\Jb)=0$ implies 
$f_i(\Phi^{-1}(\nub,\Jb))=0$ since by induction we have that $\Phi^{-1} \circ \ft_i
=f_i\circ \Phi^{-1}$ for $B'$. Hence $f_i(\Phi^{-1}(\nu,J))=\Phi^{-1}(\nut,\Jt)
=\Phi^{-1}(\ft_i(\nu,J))$, so that indeed $f_i(b)$ is defined, $f_i(b')$
and $\Phi(f_i(b))=\ft_i(\nu,J)$.

Now suppose that $f_i(b)$ is defined and $f_i(b')$ is undefined.
This implies that $b=i\otimes b'$. By induction $\ft_i(\nub,\Jb)$ is
undefined so that by Lemma~\ref{lem:varphi} we have $\overline{p}=\overline{s}$
where $\overline{p}=p_j^{(i)}(\nub)$ for large $j$ and $\overline{s}$ is the 
smallest label occurring in $(\nub,\Jb)^{(i)}$. Since $b$ is obtained from $b'$
by adding $i$ it follows that the vacancy numbers change as
$p:=p_j^{(i)}(\nu)=\overline{p}+1$ for large $j$ under $\delta^{-1}$
and the new smallest label occurring in $(\nu,J)^{(i)}$ is $s=\overline{s}$.
Hence $\widetilde{\varphi}_i(\nu,J)=p-s=1$, so that $\ft_i(\nu,J)$ is defined. 
It remains to prove that $\Phi(f_i(b))=\ft_i(\nu,J)$. Note that $f_i(b)=i+1\otimes b'$.
Let $\ell$ be the length of the largest part in $(\nub,\Jb)^{(i)}$.
Suppose that $\nub^{(i-1)}$ or $\nub^{(i+1)}$ has a part strictly bigger than $\ell$.
In this case $p_\ell^{(i)}(\nub)<\overline{p}=\overline{s}$ contradicting the fact
that $\overline{s}\le p_\ell^{(i)}(\nub)$ is the smallest label occurring in
$(\nub,\Jb)^{(i)}$. Hence both $\nub^{(i-1)}$ and $\nub^{(i+1)}$ have only parts
of length less or equal to $\ell$. Also by Lemma~\ref{lem:convex} we have
$p_\ell^{(i)}(\nub)=\overline{s}=s$ which shows that both $\delta^{-1}$ adding
$i+1$ and $\ft_i$ pick the string of length $\ell$ in $(\nub,\Jb)^{(i)}$.
Hence $\Phi(f_i(b))=\ft_i(\nu,J)$.

Let us now consider Case (2). Suppose that $\ft_i(\nu,J)$ is undefined and
$\ft_i(\nub,\Jb)$ is defined. Again by Lemma~\ref{lem:varphi} we have that
$p=s$ where $p=p_j^{(i)}(\nu)$ for large $j$ and $s$ is the smallest label
in $(\nu,J)^{(i)}$. If $\rk(\nu,J)<i+1$, then $s$ is still the smallest label in
$(\nub,\Jb)$ and by the change in vacancy numbers $\overline{p}\le p$.
Hence by Lemma~\ref{lem:varphi} $\widetilde{\varphi}_i(\nub,\Jb)=\overline{p}-s\le 0$
contradicting that $\ft_i(\nub,\Jb)$ is defined. Hence we must have
$\rk(\nu,J)\ge i+1$. In fact we want to show that $\rk(\nu,J)=i+1$.
Suppose $\rk(\nu,J)>i+1$. Then by the change in vacancy numbers by $\delta$
we have $\overline{p}=p=s$, so that $\widetilde{\varphi}_i(\nub,\Jb)=s-\overline{s}$.
So to achieve $\widetilde{\varphi}_i(\nub,\Jb)>0$ we need $\overline{s}<s$.
This can only happen if $p_{\ell^{(i)}-1}^{(i)}(\nu)=s$ and $\ell^{(i-1)}<\ell^{(i)}$.
If $m_{\ell^{(i)}-1}^{(i)}(\nu)>0$, then the string of length $\ell^{(i)}-1$ is
singular. Since $\ell^{(i-1)}<\ell^{(i)}$ this contradicts the fact that $\delta$
picks the string of length $\ell^{(i)}$ in $(\nu,J)^{(i)}$.
If $m_{\ell^{(i)}-1}^{(i)}(\nu)=0$, by convexity Lemma~\ref{lem:convexity},
we get a similar contradiction. Hence we have that $b=i+1\otimes b$.
Note that the above arguments also shows that $\widetilde{\varphi}_i(\nub,\Jb)=1$
since $\overline{s}\ge s$ and $\overline{p}=p-1$ if $\rk(\nu,J)=i+1$.
Hence $f_i(b)$ is undefined since $\varphi_i(b')=\widetilde{\varphi}_i(\nub,\Jb)
=1$.

Consider Case (2) where $f_i(b)$ is undefined and $f_i(b')$ is defined.
This implies that $b=i+1\otimes b'$. By induction $\widetilde{\varphi}_i(\nub,\Jb)
=\varphi_i(b')=1$ so that by Lemma~\ref{lem:varphi} we have $\overline{p}=\overline{s}+1$.
Hence $\widetilde{\varphi}_i(\nu,J)=p-s=\overline{p}-1-s=\overline{s}-s$ by the change 
of vacancy numbers. Therefore $\widetilde{\varphi}_i(\nu,J)=0$ if $\overline{s}=s$. It remains
to show that $p_{\ell+1}^{(i)}(\nu)\ge \overline{s}$ where $\ell:=s^{(i)}$ is the length 
of the string in $(\nub,\Jb)^{(i)}$ selected by $\delta^{-1}$.
Hence the only problem occurs if $p_{\ell+1}^{(i)}(\nub)=\overline{s}$ and
$s^{(i-1)}<\ell$. If $m_{\ell+1}^{(i)}(\nub)>0$, this means that
there is a singular string of length $\ell+1>s^{(i)}$ in $(\nub,\Jb)^{(i)}$
contradicting the maximality of $s^{(i)}$. If $m_{\ell+1}^{(i)}(\nub)=0$
one can again use convexity to arrive at similar contradiction.

By exclusion Case (3) follows from all the previous cases where at least 
one $f_i$ or $\ft_i$ is defined.

Now let $B=B^{r,1}\otimes B'$ where $r\ge 2$. Consider the following diagram:  
\begin{equation*}
\xymatrix{
 {\Path(B)} \ar[rrr]^{\Phi} \ar[ddd]_{f_i} \ar[dr] ^{\lb}& & &
        {\RC(L)} \ar[ddd]^{\ft_i} \ar[dl]_{\rclb} \\
 & {\Path(\lb(B))} \ar[r]^{\Phi} \ar[d]_{f_i} & {\RC(\lb(L))} \ar[d] ^{\ft_i}& \\
 & {\Path(\lb(B))} \ar[r]^{\Phi}  & {\RC(\lb(L))}  & \\
 {\Path(B)} \ar[rrr]_{\Phi} \ar[ur] ^{\lb}& & & {\RC(L)} \ar[ul]_{\rclb}
}
\end{equation*}
Again  the top and the bottom faces commute because of Definition~\ref{def:bij} (3).
The right face commutes by Lemma~\ref{lem:j commute}. The left face commutes
by definition of $f_i$ on the paths and we know $\lb$ is injective. By induction
hypothesis the back face commutes too. Hence the front face commutes.

Finally let $B=B^{r,s}\otimes B'$ where $s\ge 2$. Consider the following diagram:
\begin{equation*}
\xymatrix{
 {\Path(B)} \ar[rrr]^{\Phi} \ar[ddd]_{f_i} \ar[dr] ^{\ls}& & &
        {\RC(L)} \ar[ddd]^{\ft_i} \ar[dl]_{\rcls} \\
 & {\Path(\ls(B))} \ar[r]^{\Phi} \ar[d]_{f_i} & {\RC(\ls(L))} \ar[d] ^{\ft_i}& \\
 & {\Path(\ls(B))} \ar[r]^{\Phi}  & {\RC(\ls(L))}  & \\
 {\Path(B)} \ar[rrr]_{\Phi} \ar[ur] ^{\ls}& & & {\RC(L)} \ar[ul]_{\rcls}
}
\end{equation*}
As in the previous cases by Definition~\ref{def:bij} (2), Lemma~\ref{lem:i commute}
and induction hypothesis all the faces commute except the front. Since the map $\ls$
is injective the  front face of the above diagram commutes. This completes the proof
of Theorem~\ref{thm:commute}.  
\end{proof}


\begin{thebibliography}{99}

\bibitem{ASW:1999}
G.E.~Andrews, A.~Schilling, S.O.~Warnaar,   
\textit{An A$_2$ Bailey lemma and Rogers-Ramanujan-type identities},
J. Amer. Math. Soc. \textbf{12} (1999) 677--702.

\bibitem{Bailey:1949}
W.N.~Bailey,
\textit{Identities of the Rogers--Ramanujan type},
Proc. London Math. Soc. (2) \textbf{50} (1949) 1--10.

\bibitem{Bethe:1931}
H.A.~Bethe,
\textit{Zur Theorie der Metalle, I. Eigenwerte und Eigenfunktionen der
linearen Atomkette},
Z. Physik \textbf{71} (1931) 205--231.

\bibitem{Deka:2005}
L.~Deka,
\textit{Fermionic formulas for unrestricted Kostka polynomials and superconformal 
characters}, 
PhD dissertation, arXiv:math/0512536.

\bibitem{DS:2005}
L.~Deka, A.~Schilling,
\textit{New Explicit expression for $A_{n-1}^{(1)}$ supernomials},
17th International conference, FPSAC'2005, University of Messina, Italy, June 2005.

\bibitem{HKOTT:2002}
G.~Hatayama, A.~Kuniba, M.~Okado, T.~Takagi, Z.~Tsuboi,
\textit{Paths, crystals and fermionic formulae},
Prog. Math. Phys. \textbf{23} (2002) 205--272, Birkh\"auser Boston, Boston, MA.

\bibitem{HKKOTY:1999} 
G.~Hatayama, N.~Kirillov, A.~Kuniba, M.~Okado, T.~Takagi, Y.~Yamada,
\textit{Character formulae of $\hat{sl}_n$-modules and inhomogeneous paths},
Contemp. Math. \textbf{248} (1999) 243--291.

\bibitem{HKOTY:1999} 
G.~Hatayama, A.~Kuniba, M.~Okado, T.~Takagi, Y.~Yamada,
\textit{Remarks on fermionic formula},
Contemp. Math. \textbf{248} (1999) 243--291.

\bibitem{Kash:1990}
M.~Kashiwara,
\textit{Crystalizing the $q$-analogue of universal enveloping algebras},
Commun. Math. Phys. \textbf{133} (1990) 249--260.

\bibitem{KN:1994}
M.~Kashiwara, T.~Nakashima, 
\textit{Crystal graphs for representations of the $q$-analogue of classical Lie algebras},
J. Algebra \textbf{165} (1994), no. 2, 295--345.

\bibitem{Kir:2000}
A.N.~Kirillov,
\textit{New combinatorial formula for modified Hall-Littlewood polynomials},
Contemp. Math. \textbf{254} (2000) 283--333.

\bibitem{KKR:1986}
S.V.~Kerov, A.N.~Kirillov, N.Y.~Reshetikhin,
\textit{Combinatorics, the Bethe ansatz and representations of the symmetric group}
J. Soviet Math. \textbf{41} (1988), no. 2, 916--924.

\bibitem{KR:1988}
A.N.~Kirillov, N.Y.~Reshetikhin,
\textit{The Bethe ansatz and the combinatorics of Young tableaux},
(Russian) Zap. Nauchn. Sem. Leningrad. Otdel. Mat. Inst.
Steklov. (LOMI) \textbf{155} (1986), Differentsialnaya Geometriya,
Gruppy Li i Mekh. VIII, 65--115, 194;
translation in J. Soviet Math. \textbf{41} (1988), no. 2, 925--955.

\bibitem{KSS:2002}
A.N.~Kirillov, A.~Schilling, M.~Shimozono,
\textit{A bijection between Littlewood-Richardson tableaux and
rigged configurations},
Selecta Mathematica (N.S.) \textbf{8} (2002) 67--135.

\bibitem{KS:1998}
A.N.~Kirillov, M.~Shimozono,
\textit{A generalization of the Kostka-Foulkes polynomials},
J. Algebraic Combin. \textbf{15} (2002), no. 1, 27--69.

\bibitem{KOSTY:2005}
A.~Kuniba, M.~Okado, R.~Sakamoto, T.~Takagi, Y.~Yamada,
private communication.

\bibitem{LS:1978}
A.~Lascoux, M.-P.~Sch{\"u}tzenberger,
\textit{Sur une conjecture de H. O. Foulkes},
C. R. Acad. Sci. Paris S\'{e}r. A-B \textbf{286} (1978), no. 7, A323--A324.

\bibitem{MuPAD:2005}
MuPAD-Combinat available at\newline
{\tt http://www-igm.univ-mlv.fr/$\sim$descouen/MuPAD-Combinat/MuPAD-Combinat.html}.

\bibitem{NY:1997}
A.~Nakayashiki, Y.~Yamada,
\textit{Kostka polynomials and energy functions in solvable lattice models},
Selecta Math. (N.S.) \textbf{3} (1997), no. 4, 547--599.

\bibitem{Sa:2013}
R.~Sakamoto,
\textit{Rigged configurations and Kashiwara operators},
preprint arXiv:1302.4562.v1.

\bibitem{S:2005a}
A.~Schilling,
\textit{A bijection between type $D_n^{(1)}$ crystals and rigged configurations},
J. Algebra \textbf{285} (2005) 292--334.

\bibitem{S:2005}
A.~Schilling,
\textit{Crystal structure on rigged configurations},
International Mathematics Research Notices, Volume 2006, Article ID 97376, Pages 1--27.

\bibitem{SS:2001}
A.~Schilling, M.~Shimozono,
\textit{Fermionic formulas for level-restricted generalized Kostka polynomials
and coset branching functions},
Commun. Math. Phys. \textbf{220} (2001) 105--164.

\bibitem{SS:2005}
A.~Schilling, M.~Shimozono,
\textit{$X=M$ for symmetric powers},
J. Algebra \textbf{295} (2006), 562--610.

\bibitem{SW:1997}
A.~Schilling, S.O.~Warnaar,
\textit{Supernomial coefficients, polynomial identities and $q$-series},
The Ramanujan Journal \textbf{2} (1998) 459--494.

\bibitem{SW:1999}
A.~Schilling, S.O.~Warnaar,
\textit{Inhomogeneous lattice paths, generalized Kostka polynomials and A$_{n-1}$ 
supernomials},
Commun. Math. Phys. \textbf{202} (1999) 359--401. 

\bibitem{Sh:2002}
M.~Shimozono,
\textit{Affine type A crystal structure on tensor products of rectangles, Demazure 
characters, and nilpotent varieties},
J. Algebraic Combin. \textbf{15} (2002), no. 2, 151--187.

\bibitem{T:2004}
T.~Takagi,
\textit{Inverse scattering method for a soliton cellular automaton}, 
Nuclear Phys. B \textbf{707} (2005) 577--601.

\bibitem{W:2002}
S.O.~Warnaar,
\textit{The Bailey lemma and Kostka polynomials},
J. Algebraic Combin. \textbf{20} (2004) 131--171.

\end{thebibliography}
\end{document}